\documentclass[11pt,a4paper]{article}
\usepackage[T1]{fontenc}
\usepackage[utf8]{inputenc}
\usepackage[english]{babel}
\usepackage{amsmath}
\usepackage{amsfonts}
\usepackage{amssymb}
\usepackage{amsthm}
\usepackage{makeidx}
\usepackage{mathrsfs}
\usepackage{multirow}
\usepackage{color}
\usepackage[pagewise]{lineno}
\usepackage{colortbl}
\usepackage{hyperref}
\usepackage{array}
\usepackage{bigstrut}
\usepackage{setspace}
\usepackage{enumerate}
\usepackage[load=prefix,load=named,alsoload=binary]{siunitx}
\usepackage{titlesec} 
\titleformat{\section}[block]{\large\scshape\centering}{\thesection.}{1em}{} 

\usepackage[pdftex]{graphicx}
\usepackage{titling} 
\usepackage{bm}
\numberwithin{equation}{section}

\theoremstyle{definition}
\newtheorem{definition}{Definition}[section]

\theoremstyle{plain}
\newtheorem{theorem}[definition]{Theorem}

\newtheorem{proposition}{Proposition}

\theoremstyle{remark}
\newtheorem{remark}{Remark}

\numberwithin{definition}{section}
\numberwithin{proposition}{section}
\numberwithin{remark}{section}
\numberwithin{lemma}{section}

\date{}

\begin{document}  

\author{
\large{\bf{Donatella Donatelli}} \\[1ex] 
\normalsize Department of Information Engineering, Computer Science and Mathematics \\ 
\normalsize University of L'Aquila \\
\normalsize 67100 L’Aquila, Italy. \\
\normalsize \href{mailto:donatella.donatelli@univaq.it}{donatella.donatelli@univaq.it}
\and 
\large{\bf{Pierangelo Marcati}} \\[1ex] 
\normalsize GSSI - Gran Sasso Science Institute \\ 
\normalsize Viale F.\,Crispi, 7 \\
\normalsize 67100 L’Aquila, Italy.\\
\normalsize \href{mailto:pierangelo.marcati@univaq.it}{pierangelo.marcati@univaq.it}, \href{mailto:pierangelo.marcati@gssi.infn.it}{pierangelo.marcati@gssi.infn.it}
\and
\large{\bf{Licia Romagnoli}} \\[1ex] 
\normalsize Department of Information Engineering, Computer Science and Mathematics \\ 
\normalsize University of L'Aquila \\
\normalsize 67100 L’Aquila, Italy.\\
\normalsize \href{mailto:licia.romagnoli@graduate.univaq.it}{licia.romagnoli@graduate.univaq.it}
}

\title{\normalfont{A comparison of two mathematical models of the cerebrospinal fluid dynamics}} 
\maketitle
\begin{abstract}
In this paper we provide the numerical simulations of two cerebrospinal fluid dynamics models by comparing our results with the real data available in literature (see Section \ref{comparison}). The models describe different processes in the cerebrospinal fluid dynamics: the cerebrospinal flow in the ventricles of the brain and the reabsorption of the fluid.\\
In the appendix we show in detail the mathematical analysis of both models and we identify the set of initial conditions for which the solutions of the systems of equations do not exhibit blow up. We investigate step by step the accuracy of these theoretical outcomes with respect to the real cerebrospinal physiology and dynamics.\\
The plan of the paper is provided in Section \ref{plan}.
\end{abstract}

\section{INTRODUCTION}

The present paper introduces a mathematical analysis of two simplified mod\nolinebreak els describing the cerebrospinal fluid dynamics.\\
In order to clarify the framework to the reader, we shall provide the preliminary notions of anatomy and physiology on the processes underlying the craniospinal dynamics modelized by the systems of equations we are going to study. \\
Indeed, the cerebrospinal fluid (CSF) circulation network is an intricate system which surrounds the central nervous system and is incorporated into it. It has been the subject of debate since its first description in the $18^{th}$ century, emphasized by the complex vascular network of the choroid plexus that has been conventionally considered the most important structure in the production of CSF through a variety of transporters and active channels. \\
Recent outcomes in scientific methodology and outstanding improvements in the experimental tests represent fundamental tools in order to clarify the mechanisms of the CSF circulation mechanisms.
\subsection{Cerebrospinal fluid}
Cerebrospinal fluid is a clear fluid mainly composed of water ($99\%$), with the remaining $1\%$ accounted for by proteins, ions, neurotransmitters, and glucose (\cite{bulat}). It fills the surrounding spaces of the central nervous systems (CNS) of mammals and is a versatile wonder which sustains continuously the entire nervous system through the life of the organism. In the average adult man, the amount of CSF circulating at any given time is about $150$ $\ml$: in the ventricular compartments it is possible to detect the $17\%$ of the total fluid volume, the rest perfuses the cisterns and the subarachnoid space (SAS). CSF is produced at a rate of about $0,3 - 0,4$ $\ml/ \min$, i.\,e.\, $430 - 530$ $\ml$ per day (\cite{brown}). Classically the CSF flow has been thought as driven by the forces generated by heartbeats and pulmonary respiration. \\
The original theory about the CSF production establishes that the $75\%$ of all the fluid is produced by the choroid plexus epithelium, while the remaining $25\%$ arises from other CNS structures such as the ependymal wall, the cerebral parenchyma and interstitial fluid (\cite{jo}). Recently, however, there have been criticisms around the design of experiments on the choroid plexuses, questioning the thruthfulness of what we know about CSF.

\subsubsection{Cerebrospinal fluid production}
The secretion of CSF from the choroid plexus takes place as a process characterized by two different steps (\cite{brinker}). In the first stage, the plasma is passively filtered through the capillary endothelium into the choroidal interstitial space by means of the osmotic pressure gradient between the two surfaces. The ultrafiltrate fluid undergoes progressively active transport via the choroidal epithelium in the ventricular compartments.
An alternative hypothesis on the production of CSF supported by Oreskovi\'c and Klarica (\cite{oresk}) clarifies new data about the choroid plexus as the main compartment in which occurs the CSF secretion. The authors state that there are not any experiments that prove the ability of the choroid plexus of producing the expected volume of CSF. The new theory proposed describes the CSF formation as an active process that is not influenced by intracranial pressure. In balanced physiological conditions, the rate of CSF production must be equal to the rate of absorption. The authors postulate that this could extend to the flow rate, since production and absorption take place in different compartments of the circulation system. Hence, it is immediate to assert that CSF secretion is the driving force of flow and circulation if there is a constant volume of CSF.\\
The new theory moves from a more systematic approach, focusing on the perivascular spaces, which lie between the point where the cerebral vasculature descends from the SAS in the CNS through the pia mater by perforations (\cite{bul}). It is at this junction level that production and reabsorption of both interstitial and CSFs occur, due to the differences in hydrostatic and osmotic pressure between the CSF circulation system and the surrounding tissue. This means that CSF is continuously produced by means of the bloodstream and not in isolated organs involved in secretion, and any variation in CSF volume is influenced by CSF osmolarity (\cite{oresk}). \\
Although there exists evidence which supports these statements on the mixing and production of CSF, there is also extensive literature describing CSF ebbs and flows, and net flow, as well (\cite{spector}).  Indeed, according to Spector (\cite{spector}), the proposed active process in CSF formation and absorption by the entire CSF circulation system, ignores the mixing of CSF which is corroborated by the mobile cilia on the ependymal wall and by the transport of growth factors towards certain regions of the brain.

\subsubsection{CSF circulation and absorption}
\begin{figure}[h]
\centering
\includegraphics[scale=0.38]{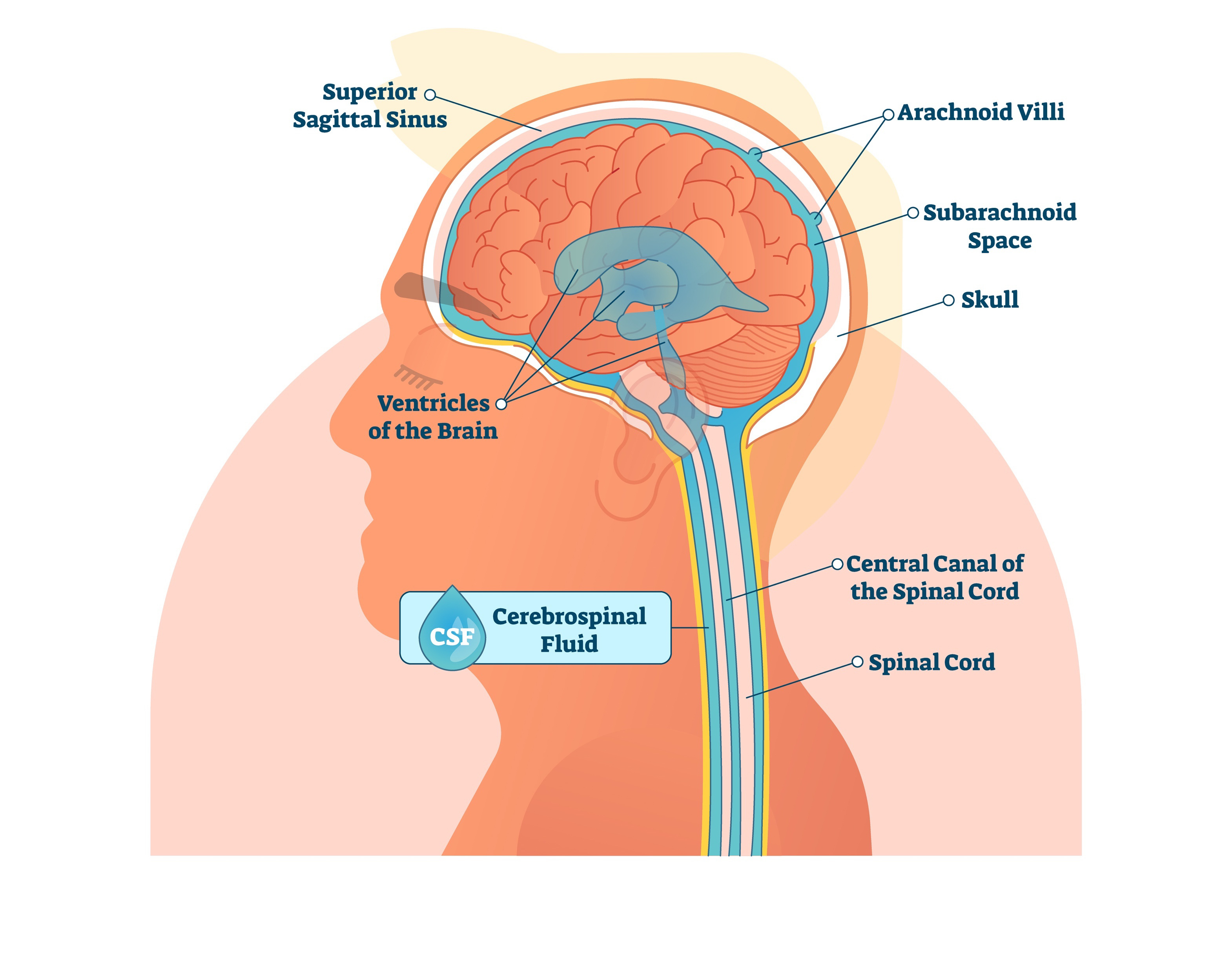}
\caption{Cerebrospinal fluid circulation. \newline{\tiny\copyright\,Copyright 2018 \url{https://www.istockphoto.com/normaals}.\label{cranio}}}
\end{figure}

After secretion, CSF generally (see Fig.\,\ref{cranio}) flows through the ventricular system, and the circulation is partially combined with the ciliated ependyma that beat in synchrony (\cite{roales}). In general it is assumed that the CSF net flow perfuses the ventricular system, by starting from the lateral ventricles (\cite{jo}).  From the lateral ventricles, the CSF flows through the left and right foramen of Monro to the third ventricle. Then, it flows downwards along the aqueduct of Sylvius in the fourth ventricle. From the fourth ventricle, the CSF is allowed to exit laterally through the foramen of Lushka, or medially through the foramen of Magendie into the SAS. The displacement through the foramen of Magendie determines the filling of the SAS. The CSF that flows along the foramen of Lushka, circulates in the subarachnoid space of the cisterns and the subarachnoid space above the cerebral cortex. Finally, the cerebrospinal fluid from the subarachnoid space is reabsorbed by arachnoid granulations which are a sort of outpouchings in the superior sagittal sinus (SSS). Arachnoid granulations act as a pathway for the reabsorption of CSF into the bloodstream through a pressure-dependent gradient (\cite{brinker}, \cite{czos}). Arachnoid granulations appear as outpouchings in the SSS due to pressure in the SAS which is greater than the pressure detected in the venous sinus (nevertheless, in vivo examination of arachnoid granulations would reveal the inverse).\\
As well as new theories on CSF formation, recently new hypothesis about CSF absorption theories have been moved (\cite{linninger2016} and references therein). Studies on rabbit and ovine models have revealed that CSF can also be significantly absorbed by means of cervical lymphatics (\cite{brinker}). CSF not reabsorbed through arachnoid granulations can reach cervical lymphatics throughout two potential pathways. The first is the SAS corresponding to the outlet of the cranial nerves (\cite{brinker}). This provides a direct path in which CSF can be shifted from cisterns to extracranial lymphatics. The second pathway across which CSF is allowed to reach the lymphatic system is along the perivascular space of the arteries and veins that penetrate the brain parenchyma (\cite{cherian}). The perivascular space (Virchow-Robin space) is the potential space which surrounds the arteries and veins of the cerebral parenchyma, that can assume different sizes depending on the pathology. When the CSF is not absorbed through the classical route, it can flow into the perivascular space or may be moved into the interstitial fluid (ISF) which is a compartment with bidirectional flow to the perivascular space and SAS. If CSF fills the ISF, finally it will be reabsorbed into the bloodstream, flows into the Virchow-Robin space or passes back to the subarachnoid space. From the perivascular space, the CSF can flow back into the subarachnoid space or be absorbed by the cervical lymphatic vessels dependent on the forces exerted by the heartbeats and pulmonary respiration.\\
Moreover, there have been conducted studies concerning the reabsorption of CSF into the dural venous plexus. At birth, the arachnoid granulations are not fully developed, and the absorption of CSF is due to the venous plexus of the inner surface of the dura which is more robust in children (\cite{mack}). Although not very wide in adult subjects, the dural venous plexus is still believed to have a role in reabsorption. Nevertheless, the exact mechanism of CSF absorption has not been completely clarified yet (\cite{mack}).

\subsubsection{CSF and cerebral vasculature}
Cerebral vasculature has the important role of providing blood supply to the brain. As the brain is one of the most important organs of the human body, an extensive network of arteries and veins assolves the task of providing oxygenated blood to the organ and to drain deoxygenated blood from it. The arterial system supplying the brain presents many collateral vessel that insert redundancies into the system. \\
The cerebral venous compartment is believed to be the main responsible of the intracranial compliance (defined as the ratio of volume and pressure change). In fact, the possibility of emptying part of the venous compartment to make up for an increased volume of the other compartment is a key phenomenon in the dynamic equilibrium of intracranial pressure (ICP). Moreover, an increased ICP can cause an unstable collapse of veins and so induce a change of resistance that may alter the cerebral blood flow (CBF). Such modification of the resistance may trigger autoregulation mechanisms, which act to ensure a constant blood flow, and it may be responsible for the amplification of the CSF pulsewave observed during infusion studies. However, no definite conclusions have been drawn on the latter point.

One of the main feature of the cerebral vasculature involves the coupling with the cerebrospinal fluid circulation.\\
The CSF bulk flow associated to the production and reabsorption is small compared to its pulsatile component. The oscillations of pulsatile CSF are assumed to be managed by systolic vasculature dilatation followed by diastolic contraction. \\
During the normal cardiac cycle, approximately $750$ $\ml$ of blood are pushed into the head every minute. The increase in systolic blood pressure inflates arterial blood vessels, therefore an increase in the volume of brain blood during systole is detected. Since the cranial vault is characterized by a rigid skull in adults, the vascular expansion of the major arteries which pass through the spaces perfused by the cerebrospinal fluid activates the CSF motion. MRI techniques show evidences that the total volume of cerebral blood expands and comprises itself in every cardiac cycle of about $1-2$ $\ml$, the same volumetric quantity corresponding to the exchange of CSF between the cranial and spinal SAS (\cite{linninger2016} and references therein).\\
The expansion of the vascular volume could be transmitted from the cortical surface through brain tissue, whose compression causes the contraction of the ventricular space. Otherwise, the movement of the ventricular wall could derive from inside the ventricles by systolic expansion of the choroidal arteries, such that the ventricular walls pulsate against the periventricular ependymal layer. Ventricular dilatation due to the expansion of the choroid was proposed in the theoretical model of Linninger et al.\,(see \cite{b3}).

Alternatively, Buishas et al.\,(\cite{bua}) proposed a model of water transport through the parenchyma from the microcirculation as driven by Starling forces. This model investigates the effect of osmotic pressure on water
transport between the cerebral vasculature, the extracellular space, the perivascular space and the CSF and predicts the effects the osmolarity of ECS, blood, and CSF on water flux in the brain, establishing a link between osmotic imbalances and pathological conditions such as hydrocephalus and edema.

\subsection{Key points of the paper}
In the light of the above, the main issues in analyzing the intracranial dynamics is represented by the particular structure of the human brain that includes very complex elements and the occurrence of many different phenomena:
\begin{itemize}
\item The flow of the CSF fluid throughout the CSF compartments. 
\item The mechanical interaction between the fluid and the brain. The brain parenchyma and the CSF compartments
interact by exerting reciprocal stress on each other. 
\item Coupling with the circulatory-system. Models of the global circulatory tree may be necessary to provide pressure values in the brain with enough accuracy and to impose the Monro-Kellie doctrine on arterial and venous system, CSF and brain parenchyma.\\
The Monro-Kellie doctrine, may be expressed mathematically through the Linninger reformulation (see \cite{linninger2009}):
\begin{align}\label{monro}
V&=\sum_b V_b+\sum_{CSF}+V_{brain}=\mbox{constant}, \notag\\
V_{brain}&=V_{exf}+V_{solidBrain}
\end{align}
The subscript, exf, refers to the extracellular fluid flow inside brain. Moreover, the volume of the brain parenchyma in $\eqref{monro}_2$ is the sum of the volume of extracellular fluid, $V_{exf}$, and the volume of the solid part of the brain, $V_{solidBrain}$. The first quantity in $\eqref{monro}_1$ is characterized by the sum of all the volumes of the vascular system components (i.\,e.\,arteries, arterioles, capillaries, veinules, veins and venous sinus), while the second quantity is the sum of the volumes related to the compartments involved in the ventricular system (i.\,e.\,lateral ventricles, third and fourth ventricles, cranial subarachnoid space).
\item Production and reabsorption laws. In order to describe the CSF system, suitable laws are needed to model the production of fluid inside the ventricles and the absorption of the fluid in the sagittal sinus.
\end{itemize}

Lumped or compartmental models are particularly appreciated in the description of complex multicompartmental systems, in fact they are usually simple to develop as prototype and fast to solve.\\
Marmarou pioneered the field of mathematical modeling of intracranial dynamics by lumped models in 1973. In his thesis \cite{marmarou}, he introduced the concept of PVI (pressure-volume index), studied the intracranial pressure-volume relation, defined how to assess PVI and CSF reabsorption resistance by infusion studies and developed the first lumped model for intracranial dynamics focused on CSF \cite{marmarou2}. His aim was to predict changes in ICP when the content of CSF inside the CSF system changes.\\

Brain parenchyma has been modelled through three dimensional poroelastic models (\cite{chou2016}, \cite{guo2018}, \cite{chou2014}). Brain fluid systems have been modelled through lumped parameter models (\cite{ursino1988}, \cite{ursino88}, \cite{gadda2015}, \cite{gehlen2017}) and through multi-scale models (\cite{muller2014}). Ursino (\cite{ursino1988}, \cite{ursino88}) proposed a mathematical model of the human intracranial hydrodynamics. The group of Linninger (\cite{b3}) created a dynamic model of the ventricular system to test the hypothesis that choroid plexus expansion drives CSF flow in this system. Later, they proposed a mathematical model of blood, cerebrospinal fluid and brain dynamics, including the Monro-Kellie doctrine (\cite{linninger2009}). The same group presented a mathematical model of the intracranial fluid dynamics based on the Bulat-Klarica-Ore\v{s}kovi\'c hypothesis 
(\cite{oreskovic2017}, \cite{linninger2017}). Gehlen et al.\,(\cite{gehlen2017}) studied the effect of postural changes in the CSF dynamics through a lumped-parameter model of the CSF system and major compartments of the cardiovascular system.\\

Our aim is to compare from a mathematical point of view two mod\nolinebreak els which describe the simplified dynamics of the cerebrospinal fluid in the parenchyma. Both models are obtained by considering as a starting point the CSF model studied in \cite{donmaro}, where the authors provide a comprehensive mathematical analysis of the system of equations derived from the CSF compartmental model introduced by Linninger et al.\,in \cite{b3}. \\
\begin{figure}[t]
\centering
\includegraphics[scale=0.47]{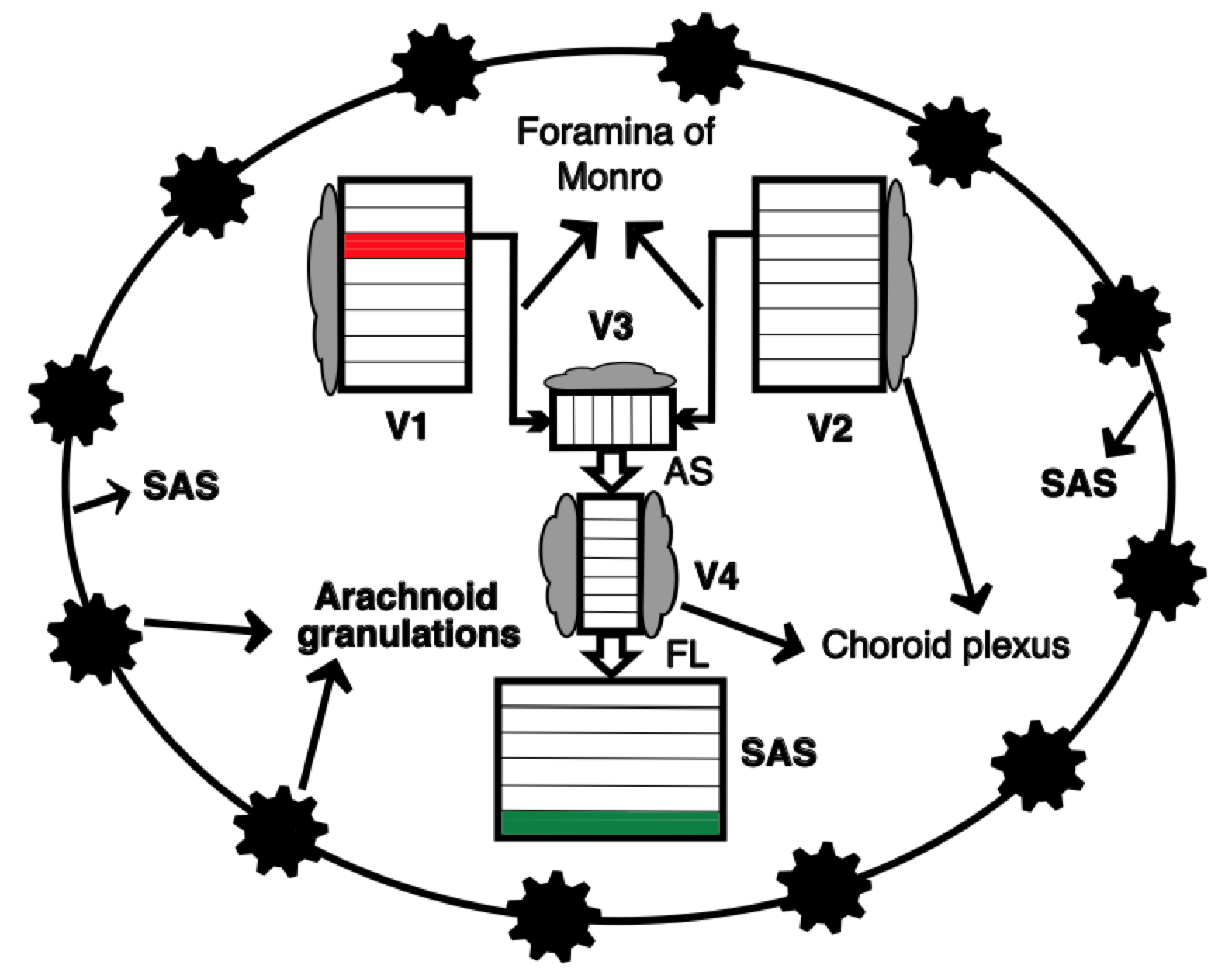}
\caption{Discretization of CSF compartments: CSF secreted by the choroid plexi (in grey) flows into the lateral ventricles (\textbf{V1} e \textbf{V2}), then, through the foramina of Monro, it fills the third ventricle (\textbf{V3}) and flows downwards into the fourth ventricle (\textbf{V4}) through the Acqueduct of Sylvius (\textbf{AS}). It passes into the Foramen of Luschka (\textbf{FL}) and perfuses the subarachnoidal space (\textbf{SAS}) where it is partially reabsorbed by the arachnoid granulations.\label{discretizzato}}
\end{figure}
The decision to adopt non recent CSF models (\cite{b3}, \cite{marmarou2}) as starting point for our research, was dictated by the fact that the analysis that we propose below represents a first attempt to work from a purely mathematical point of view with models that describe cerebrospinal fluid dynamics. 
Because of the complexity of the mechanisms involved in CSF production, circulation and absorption, we decide to neglect in this first step some fundamental contributions of the brain vasculature or the parenchyma. Therefore, the aim of this paper is to lay the foundations for a much more detailed future study, based on models that treat in detail the interactions that emerge in the intracranial pattern among the vascular system, the CSF compartments and the brain parenchyma, and are able to predict blood alterations or potential cerebral water content variation. \\
In this paper we want to improve the simplified model analyzed in \cite{donmaro} in the following way:
\begin{itemize}
\item the cross sectional area of the CSF compartments, which is affected by the CSF flow and by the intracranial pressure in the real physiology, is not assumed anymore as a constant but it is free to change during the cerebral process modelized;
\item we add to the Linninger model (see \cite{b3}) an equation, developed by Marmarou et al. in \cite{marmarou2}, which describes the evolution of the intracranial pressure;
\item we introduce a CSF absorption rate in order to study also the equations that rule the final part of the CSF process which involves the venous superior sagittal sinus. We will define this new quantity in terms of intracranial pressure. 
\end{itemize}
A model which reflects a perfect mixture of these features is the one proposed by Linninger et al.\, in 2009 (see \cite{linninger2009}) that provides a complete representation of the CSF pathways interacting with the vascular system and the porous parenchyma. 

\subsection{Definition of equations symbols}\label{simboli}

\begin{figure}[htbp]
\centering
\includegraphics[scale=0.75]{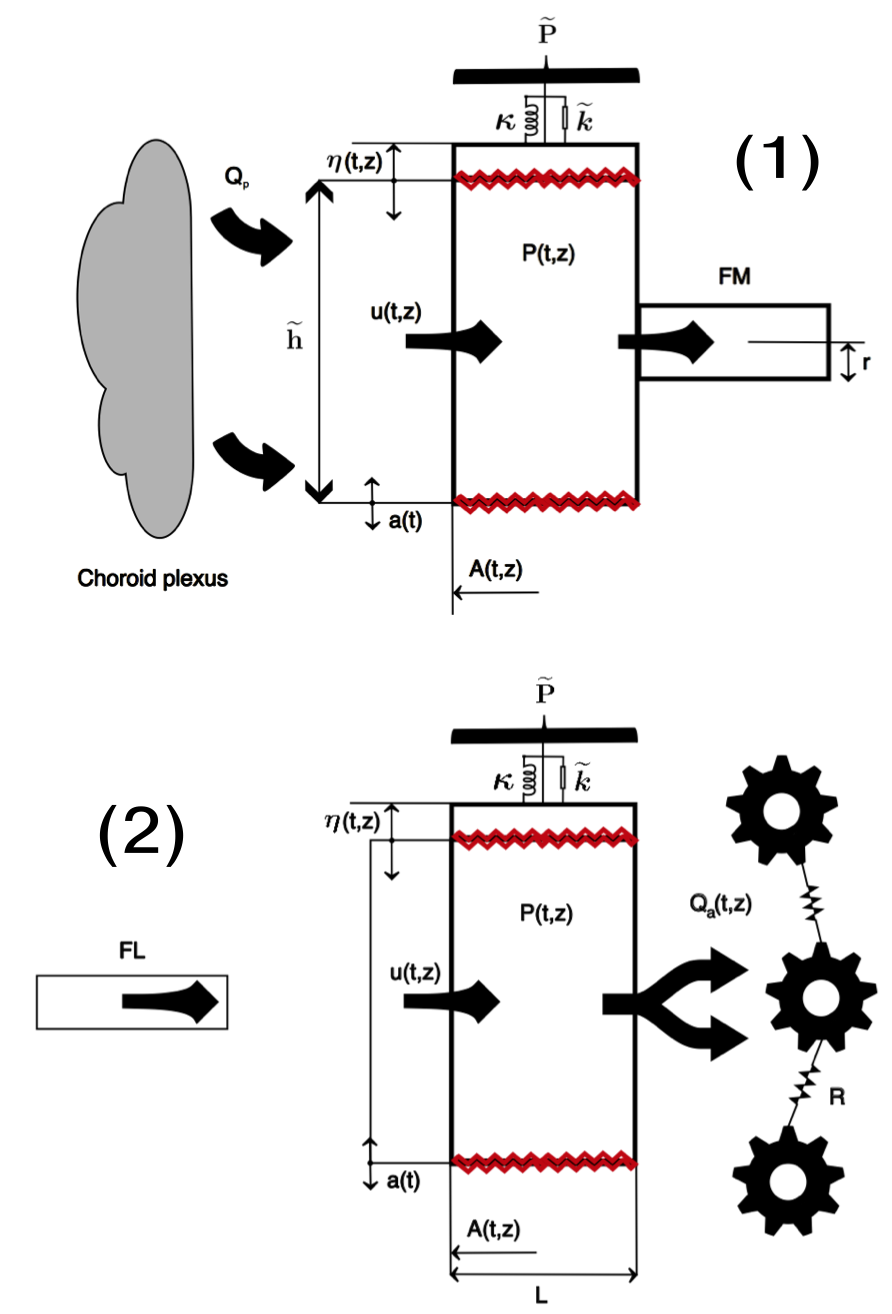}
\caption{$\mathbf{(1)}$ Zoom on the red cylinder of the lateral ventricle in Fig.\,\ref{discretizzato}. $\mathbf{(2)}$ Zoom on the green cylinder of the SAS in Fig.\,\ref{discretizzato}. Figure adapted from Linninger et al.\,(\cite{b3}).\label{zooms}}
\end{figure}

The quantities involved in the models are the following (see Fig.\,\ref{zooms}).

\begin{table}[h]
\centering
\begin{tabular}{|r|c|c|}
\hline
&\textbf{Definition}&\\ \hline
$A(t,z)$&\textbf{cross section} of the CSF compartments & $\left[\meter^2\right]$\\ \hline
$\eta(t,z)$&\textbf{tissue displacement} in a section&$[\meter]$\\ \hline
$Q_a(t,z)$&\textbf{CSF absorption rate}& $\left[\ml/\min\right]$ \\ \hline
$u(t,z)$& axial CSF \textbf{flow velocity}& $\left[\meter/\second\right]$ \\ \hline
$P(t,z)$& \textbf{CSF pressure} in ventricles and SAS (ICP)&$\left[\N/\meter^2\right]$\\ \hline
\end{tabular}
\caption{CSF dynamics quantities involved in Models A1 and A2.}
\end{table}

The \textbf{time} $t$ is such that $t\in[0,T_0]$ and $z\in[0,L]$ is the \textbf{axial coordinate}.\\
The physical constants that will be adopted in the present paper are listed in Table \ref{costanti} (see also Fig.\,\ref{zooms}).\\
The forcing function given by 
\begin{displaymath}
a(t)=\overline{\alpha}\left(1.3+\sin\left(\omega t-\frac{\pi}{2}\right)-\frac{1}{2}\cos\left(2\omega t-\frac{\pi}{2}\right)\right),
\end{displaymath}
represents the choroid plexus periodic motion which follows the cardiac cycle.
\vspace{4pt}
\begin{table}[h]
\centering
\begin{tabular}{|r|c|c|}
\hline
&\textbf{Definition}&\\ \hline
$\rho$&fluid density&$\left[\kg/\meter^3\right]$\\ \hline
$\delta$&tissue width&$\left[\meter\right]$\\ \hline
$\widehat{\alpha}=\rho\delta$&&$\left[\kg/\meter^2\right]$\\ \hline
$\kappa$&tissue elasticity constant&$\left[\N/\meter\right]$\\ \hline
$\widetilde{k}$&tissue compliance&$\left[\N\,\second/\meter\right]$\\ \hline
$r$&radius of the foramina and aqueduct&$[\meter]$\\ \hline
$\mu$&fluid viscosity&$\left[\Pa\,\second\right]$\\ \hline
$\beta$&Poiseuille friction term ($8\mu/r^2$)&$\left[\N/\meter^3\right]$\\ \hline
$\widetilde{h}$&height of the CSF compartments&$[\meter]$\\ \hline
$Q_{p}$&CSF production rate in the choroid plexus&$\left[\meter^3/\second\right]$\\ \hline
$R$&resistance to CSF absorption&$\left[\Pa\,\second/\meter^3\right]$\\ \hline
$\overline{\alpha}$&amplitude of choroid expansion&$[\meter]$\\ \hline
$\omega$&heart rate frequency&$\left[\rad/\second\right]$\\ \hline
$K$&mathematical constant&$\ml$\\ \hline
$\widetilde{P}$&pressure of brain parenchyma&$\left[\N/\meter^2\right]$\\ \hline
$L$&length of a single discretized cylinder&$[\meter]$\\ \hline
$a(t)$&forcing function&$[\meter]$\\ \hline
\end{tabular}
\caption{Physical constants.\label{costanti}}
\end{table}

\subsection{Description of the models}

This work focuses on two systems of equation given by a CSF compartmental model (see \cite{b3}) where the CSF sections are described by a cylindrical discretization with axial symmetry and radial displacements. The foramina are treated as elastic pipes. Since the CSF flow is basically laminar based on its very low Reynolds numbers $(Re<100)$, the flow friction was expressed as a function of the cross-sectional area. The latter, denoted by $A(t,z)$, is considered as a function depending on time and space, and this represents the first significant improvement on the analysis developed in \cite{donmaro}, where it is assumed to be constant.\\
As already explained, the models we will analyze in the present paper derive from the compartmental model presented by Linninger et al.\, in 2005 (\cite{b3}) which the authors improved in 2009 (\cite{linninger2009}), by proposing a comprehensive model of human intracranial dynamics where cerebral blood, CSF and brain parenchyma as well as the spinal canal were included. \\
In order to manage properly our mathematical analysis, in this first approach to the CSF dynamics, we take into account the first model mentioned since it represents a good starting point for our analysis at this early stage.\\
For the same reason, to define the behaviour of the intracranial pressure, we consider the model proposed by Marmarou et.\,al in \cite{marmarou2}, described by the following differential equation that rules the CSF hydrodynamic without incorporating the interactions with brain vasculature and porous parenchyma,
\begin{displaymath} 
\partial_t P(t,z)-\displaystyle{\frac{K}{R}}P(t,z)^2-KP(t,z)\left(Q_p+\displaystyle{\frac{\widetilde{P}}{R}}\right)=0.
\end{displaymath}
Therefore, the first simplified CSF model we study in this paper reads as follows,
\begin{flushleft}
\textbf{Model A1}
\end{flushleft}
\begin{eqnarray}\label{br1}
\begin{cases}
\widetilde{h}\partial_tA(t,z)+\partial_t\left(A(t,z)a(t)\right)+\partial_t(A(t,z)\eta(t,z))+A(t,z)u(t,z)-Q_p=0, \\
\widehat{\alpha} A(t,z)\partial_{tt}\eta(t,z)+\widetilde{k}\partial_t\eta(t,z)+\kappa\eta(t,z)-A(t,z)P(t,z)+A(t,z)\widetilde{P}=0,  \\
\rho\partial_t u(t,z)+\rho u(t,z)\partial_z u(t,z)+\beta u(t,z)+\partial_z P(t,z)=0,  \\
R\partial_t P(t,z)-KP(t,z)^2-KP(t,z)\left(RQ_p+\widetilde{P}\right)=0.
\end{cases}
\end{eqnarray}

The previous model is designed as follows.
\begin{itemize}
\item The first equation of \textbf{Model A1} is the continuity equation, consistent with the assumption of incompressible CSF. It defines the conservation of fluid mass through the motion. 
\item The second equation represents the distensibility balance, indeed the stresses, strains and displacements detectable in parenchyma and in the CSF compartments tissue, is described by applying the elastodynamics laws (see \cite{b3}).
\item The third equation is the axial momentum balance where $\beta u(t,z)$ is the Poiseuille term (see \cite{b3}).
\item The last equation is the intracranial pressure equation already defined (\cite{marmarou2}).
All the variables and constants that appear in the equations are defined in detail in Section \ref{simboli}.\\
\end{itemize}
As previously explained, we know that under normal physiological conditions, the rate of CSF formation ($Q_p$) is balanced by an equal rate of absorption ($Q_a$). This condition of equilibrium results in no increase or decrease in the amount of volume stored and the initial resting volume as well as the CSF pressure ($P(t,z)$) are maintained at a constant level. The rate of outflow ($Q_a$) is given by the gradient of pressure between CSF space and the venous system of the dural sinus (we assume equal to $\widetilde{P}$) divided by the resistance to absorption ($R$) (see \cite{marmarou2}):
\begin{displaymath}
Q_a(t,z)=\displaystyle{\frac{1}{R}}\left(P(t,z)-\widetilde{P}\right).
\end{displaymath}
This last issue allows us to introduce the second model we want to study and compare with \eqref{br1}.

\begin{flushleft}
\textbf{Model A2}
\end{flushleft}
\begin{eqnarray}\label{br2}
\begin{cases}
\widetilde{h}\partial_tA(t,z)+\partial_t\left(A(t,z)a(t)\right)+\partial_t(A(t,z)\eta(t,z))+Q_a(t,z)-Q_p=0, \\
RQ_a(t,z)-\left(P(t,z)-\widetilde{P}\right)=0,\\
\widehat{\alpha} A(t,z)\partial_{tt}\eta(t,z)+\widetilde{k}\partial_t\eta(t,z)+\kappa\eta(t,z)-A(t,z)P(t,z)+A(t,z)\widetilde{P}=0,  \\
\rho\partial_t u(t,z)+\rho u(t,z)\partial_z u(t,z)+\beta u(t,z)+\partial_z P(t,z)=0,  \\
R\partial_t P(t,z)-KP(t,z)^2-KP(t,z)\left(RQ_p+\widetilde{P}\right)=0.
\end{cases}
\end{eqnarray}

The previous model, along the lines of \textbf{Model A1} is defined as follows.
\begin{itemize}
\item The first equation is the continuity equation where we are taking into account a simplified reabsorption term, $Q_a$.
\item The second equation defines the reabsorption term (\cite{marmarou2}).
\item The third equation is the elastodynamics equation (\cite{b3}).
\item The fourth equation is the momentum balance.
\item The last equation is the intracranial pressure equation (\cite{marmarou2}).
\end{itemize}

In order to simplify the models, we reduce the order of the second equation of \eqref{br1} and \eqref{br2}. Furthermore, in the system \eqref{br2} we plug the equation 
\begin{displaymath}
Q_a(t,z)=\displaystyle{\frac{1}{R}}\left(P(t,z)-\widetilde{P}\right)
\end{displaymath}
in the equation $\eqref{br2}_1$. \\
Therefore we reformulate the models in the following way 
\begin{flushleft}
\textbf{Model A1}
\end{flushleft}
\begin{eqnarray}\label{br1bis}
\begin{cases}
\partial_t\eta(t,z)=\zeta(t,z),\\
\widehat{\alpha} A(t,z)\partial_{t}\zeta(t,z)+\widetilde{k}\zeta(t,z)+\kappa\eta(t,z)-A(t,z)P(t,z)+A(t,z)\widetilde{P}=0, \\
\left(\widetilde{h}+a(t)+\eta(t,z)\right)\partial_t A(t,z)+\left(a'(t)+\zeta(t,z)+u(t,z)\right)A(t,z)-Q_p=0,\\
\rho\partial_t u(t,z)+\rho u(t,z)\partial_z u(t,z)+\beta u(t,z)+\partial_z P(t,z)=0,  \\
R\partial_t P(t,z)-KP(t,z)^2-KP(t,z)\left(RQ_p+\widetilde{P}\right)=0,
\end{cases}
\end{eqnarray}
\newpage
\begin{flushleft}
\textbf{Model A2}
\end{flushleft}
\begin{eqnarray}\label{br2bis}
\begin{cases}
\partial_t\eta(t,z)=\zeta(t,z),\\
\widehat{\alpha}A(t,z)\partial_{t}\zeta(t,z)+\widetilde{k}\zeta(t,z)+\kappa\eta(t,z)-A(t,z)P(t,z)+A(t,z)\widetilde{P}=0, \\
\left(\widetilde{h}+a(t)+\eta(t,z)\right)\partial_t A(t,z)+\left(a'(t)+\zeta(t,z)\right)A(t,z)-Q_p\\
\ \ \ \ \ \ \ \ \ \ \ \ \ \ \ \ \ \ \ \ \ \ \ \ \ \ \ \ \ \ \ \ \ \ \ \ \ \ \ \ \ \ \ \ \ \ \ \ \ \ \ \ \ \ \ \ \ \ \ \ \ \ \ \ \ \ \ +\displaystyle{\frac{1}{R}}\left(P(t,z)-\widetilde{P}\right)=0, \\
\rho\partial_t u(t,z)+\rho u(t,z)\partial_z u(t,z)+\beta u(t,z)+\partial_z P(t,z)=0,  \\
R\partial_t P(t,z)-KP(t,z)^2-KP(t,z)\left(RQ_p+\widetilde{P}\right)=0,
\end{cases}
\end{eqnarray}
where $a'(t)=\displaystyle{\frac{d}{dt}}a(t)$.\\
In both models we impose for the unknown $u(t,z)$ the following initial and boundary conditions respectively
\begin{equation}\label{br3}
u(0,z)=u_0(z)=f(z)\,\in\,\mathcal{H}^{s}([0, L]),
\end{equation}
where $\displaystyle{s>\frac{5}{2}}$ and
\begin{equation}\label{br4}
u(t,0)=u(t,L)=0,
\end{equation}
and the following compatibility conditions
\begin{align}\label{br5}
u(0,0)&=u_0(0)=f(0), \notag \\
u(0,L)&=u_0(L)=f(L),
\end{align}
for which we assume
\begin{equation}\label{br6}
f(0)=f(L)=0.
\end{equation}
Similarly for the unknown $\eta$ we assume
\begin{equation}\label{br7}
\eta(0,z)=\eta_0(z)=g(z)\in\mathcal{H}^{s}([0,L])
\end{equation}
and the following compatibility conditions
\begin{align}\label{br8}
\eta(0,0)&=\eta_0(0)=g(0),  \notag \\
\eta(0,L)&=\eta_0(L)=g(L).
\end{align}
Concerning the no-slip conditions, we can assume, without loss of generality, the following boundary conditions
\begin{equation}\label{br11}
\eta(t,0)=\eta(t,L)=0,
\end{equation}
then we get
\begin{equation}\label{br12}
g(0)=g(L)=0.
\end{equation}
Moreover for the extra unknown $\zeta$ we impose
\begin{equation}\label{br13}
\zeta(0,z)=\zeta_0(z)=q(z)\in\mathcal{H}^{s}([0,L])
\end{equation}
as initial condition, and the following compatibility conditions
\begin{align}\label{br14}
\zeta(0,0)&=\zeta_0(0)=q(0),  \notag \\
\zeta(0,L)&=\zeta_0(L)=q(L).
\end{align}
By considering the boundary conditions \eqref{br11} of the unknown $\eta$, we get 
\begin{equation}\label{br15}
\zeta(t,0)=\zeta(t,L)=0
\end{equation}
and
\begin{equation}\label{br16}
q(0)=q(L)=0.
\end{equation}

In both systems we have a homogeneous Riccati equation that describes the evolution of the intracranial pressure. For the latter we choose the following initial condition
\begin{equation}\label{br17}
P(0,z)=P_0(z)=b(z)\in\mathcal{H}^{s}([0,L]).
\end{equation}
By using classical ODE methods we can find the general solution of the Riccati equation 
\begin{equation}\label{br18}
P(t,z)=\displaystyle{\frac{b(z)e^{\mathcal{C}t}}{1+\displaystyle{\frac{K}{R\mathcal{C}}}b(z)\left(1-e^{\mathcal{C}t}\right)}},
\end{equation}
where $\mathcal{C}=\displaystyle{K\left(Q_p+\frac{\widetilde{P}}{R}\right)}$ and $t\in[0,T_0]$.\\
It represents the variation in time of the pressure during the entire CSF process from the choroid plexus to the subarachnoid space.
\begin{remark}\label{remark1}
It is important to observe that the measurement of the CSF pressure gradients can be affected by the body position. Since the common techniques are rather invasive, we have only limited data related to the upright posture. Indeed, CSF pressure is usually measured while a person is lying in a horizontal supine position. Normal CSF pressure values, in that case, are around $7 - 15$ $\mmHg$, and the pressure is the same along the SAS and inside the cranial vault (see \cite{davson}). It is known that the change in body position (from horizontal to upright, head up or sitting position) is followed by a drop in intracranial pressure to the subatmospheric value, and it can be detected a new pressure gradient along the craniospinal axis. In \cite{alperin} and \cite{klarica}, the authors carryed out analysis using data previously collected from healthy subjects scanned in supine and sitting positions and, after modifications for the hydrostatic component, they detected ICP values at the central cranial location of all subjects were negative, with an average value of $-3,4$ $\mmHg$.\\
By taking into account the previous considerations, it is clear that there are two different choices of the initial condition for the intracranial pressure and they depend on which situation we want to investigate. \\
More precisely:
\begin{itemize}
\item if we apply Models A1 and A2 to a subject that is lying supine, we are forced to fix an initial datum $b(z)>0$. In this case, the solution \eqref{br18} blows up at the finite time 
\begin{equation}\label{br19}
\widetilde{T}_0=\displaystyle{\frac{1}{\mathcal{C}}\ln\left(\frac{R\mathcal{C}}{Kb(z)}+1\right)},
\end{equation}
and, in order to study local existence in time of the solutions of the systems \eqref{br1bis} and \eqref{br2bis}, henceforth we will consider the time $t\in[0,T_0]$ where $T_0<\widetilde{T}_0$.\\
We point out that this is not a big restriction on the possible values of $t$, since if we plug the correct values of the physical parameters in \eqref{br19} and we take
\begin{displaymath}
b(z)=7\,\mbox{mmHg}= 932,54\,\mbox{Pa},
\end{displaymath}
which is the average intracranial pressure, we can observe that the blow up time $\widetilde{T}_0$ corresponds to an age of 124 years.
\item If we consider a subject in the upright position, we need to choose the initial datum $b(z)<0$. This implies there is no blow up at finite time and no restriction on the time interval is required. 
\end{itemize}
\end{remark}
The last unknown to analyze is the generic axial section which is denoted by $A(t,z)$. We impose the following initial condition
\begin{equation}\label{br23}
A(0,z)=A_0(z)=h(z)\in\mathcal{H}^{s}([0,L]),
\end{equation}
and we observe that the general solution of the third equation in both systems \eqref{br1bis} and \eqref{br2bis} is
\begin{equation}\label{br26}
A(t,z)=h(z)e^{G(t,z)t}+e^{-G(t,z)t}\int_0^t H(s,z)e^{G(s,z)s}\,ds,
\end{equation}
where 
\begin{equation}\label{br27}
G(t,z)=\displaystyle{\frac{\displaystyle{a'(t)}+\zeta(t,z)+u(t,z)}{\widetilde{h}+a(t)+\eta(t,z)}}, \ \ \ \ 
H(t,z)=\displaystyle{\frac{Q_p}{\widetilde{h}+a(t)+\eta(t,z)}},
\end{equation}
for Model A1, and
\begin{equation}\label{br28}
G(t,z)=\displaystyle{\frac{\displaystyle{a'(t)}+\zeta(t,z)}{\widetilde{h}+a(t)+\eta(t,z)}}, \ \ \ \ H(t,z)=\displaystyle{\frac{RQ_p-P(t,z)-\widetilde{P}}{R(\widetilde{h}+a(t)+\eta(t,z))}},
\end{equation}
for Model A2. 

\vspace{12pt}

Our approach in the analysis of Models A1 and A2, supplemented by the initial and boundary conditions previously defined is as follows.
\begin{itemize}
\item We implement a procedure of successive approximations to the systems A1 and A2 and we analyze the regularity of the sequences $u^n(t,z),$ $P^n(t,z), \eta^n(t,z), \zeta^n(t,z)$ and $A^n(t,z)$ in order to prove the existence of an approximating solution. The sequence $u^n$ requires a different analysis due to the nonlinearity of the velocity flow equation.
\item We show the convergence of the sequence $u^n$ for which we resort to high order energy estimates. For the other sequences, in order to obtain convergence to classic solutions, we need to apply first the Ascoli's Theorem \ref{ascoli} which guarantees that the set 
\begin{displaymath}
\Psi(t,z)=\left\{P^n(t,z), \eta^n(t,z), \zeta^n(t,z), A^n(t,z); t\in[0,T], z\in[0,L]\right\}
\end{displaymath}
is relatively compact in $C([0,T]; \mathcal{H}^{s-1}([0,L]))\cap C^1((0,T]; \mathcal{H}^{s-2}([0,L]))$ provided the sequences are equicontinuous and equibounded. 
\item We pass into the limit in the approximating systems and we complete the proof of the local existence and uniqueness of classical solutions for Models A1 and A2.
\item We investigate under which conditions it is possible to obtain global existence. Therefore we study the equation $\eqref{br1bis}_4$ with the standard characteristic method and we obtain a nonhomogeneous Riccati equation (NRE) concerning which we have to find a general solution. In constructing the particular real solution to the NRE, we are forced to employ a condition on the pressure; this result allows us to prove global existence in time and uniqueness of solutions for the systems A1 and A2 under precise restrictions on the initial data.
\item In order to prove that the analysis of Model A1 and A2 is reliable, we implement proper numerical simulations by fixing first initial data which fulfill the condition imposed for the global existence of the solutions and then we choose initial data which violate it.
\end{itemize}

\subsection{Plan of the paper} \label{plan}
The plan of the present manuscript is as follows. In Section 2 we collect all the definitions and the technical results we are going to use through the paper and we state our main results which are the local existence Theorems \ref{TeoremaBr1} and \ref{TeoremaBr2} and the global existence Theorem \ref{TeoremaBr3}. In Section 3, we perform the numerical simulations that validate the results on the global existence of the solutions obtained for Models A1 and A2 and we verify for each simulation how the behavior of the quantities involved in both models matches the real CSF dynamics. Section 4 is entirely devoted to the comparison of our theoretical results with respect to the real data available in literature. Finally, in the Appendix we provide the detailed mathematical analysis of both models in order to prove the local existence Theorems \ref{TeoremaBr1} and \ref{TeoremaBr2} and the global existence of solutions (Theorem \ref{TeoremaBr3}) for which we are able to show the set of initial data that avoids blow up formation. 

\section{PRELIMINARIES}

In this section we are going to fix the notations used in the paper, we recall the main tools we need to study our problem and we state our main result.

\subsection{Notations and definitions}\label{notazioni}

In the sequel we shall use the customary Lebesgue spaces $L^p(\Omega)$ and Sobolev spaces $\mathcal{H}^s(\Omega):=W^{s,2}$ with $\Omega:=[0,L]$. We use $\left\|\cdot\right\|_{L^{\infty}_t\mathcal{H}^s_z}$ to denote the norm $\left\|\cdot\right\|_{L^{\infty}((0,T))\mathcal{H}^{s}([0,L])}$ and $\left\|\cdot\right\|_2$ to denote the norm $\left\|\cdot\right\|_{L^2([0,L])}$.\vspace{4pt}
\\
We will use also the following Sobolev interpolation theorem
\begin{equation}\label{ar2}
\left\| w\right\|_{\mathcal{H}^{r'}([0,L])}\leq C_r\left\| w\right\|_{L^{2}([0,L])}^{1-\frac{r'}{r}}\left\| w\right\|_{\mathcal{H}^{r}([0,L])}^{\frac{r'}{r}},
\end{equation}
which holds for any $0\leq r'\leq r$ and $w\in\mathcal{H}^r([0,L])$.\\

Now, we recall a theorem and some definitions (\cite{stampacchia}) which will be crucial for the proof of the local existence Theorems \eqref{TeoremaBr1} and \eqref{TeoremaBr2}.
\begin{theorem}\label{iperb}
Let $\mathbb{A}_j(x,t),\, j=0,..., N$ be symmetric $m\times m$ matrices in $M\times (0, T)$ with $M\subset\mathbb{R}^N$, $T>0$, $\mathbf{f}(x, t)$ and $\mathbf{v}_0(x)$, m-dimensional vector functions defined in $M\times (0, T)$ and $M$, respectively. Let us consider the problem 
\begin{subequations}
\begin{equation}\label{ar44}
\mathbb{A}_0(x,t)\partial_t\mathbf{v}+\sum_{j=1}^N\,\mathbb{A}_j(x,t)\partial_j\mathbf{v}+\mathbb{B}(x,t)\mathbf{v}=\mathbf{f}(x,t), 
\end{equation}
\begin{equation}\label{ar45}
\mathbf{v}(x,0)=\mathbf{v}_0(x).
\end{equation}
\end{subequations}
Let
\begin{displaymath}
\mathbb{A}_j\in\left[C([0,T], \mathcal{H}^s(M))\cap C^1((0,T), \mathcal{H}^{s-1}(M))\right]^{m\times m}, \ \ j= 0, 1,..., N,
\end{displaymath}
$\mathbb{A}_0(x,t)$ invertible, $\inf_{x,t}\left\|\mathbb{A}_0(x,t)\right\|_{M\times M}>0$, 
$\mathbb{B}\in C((0,T); \mathcal{H}^{s-1}(M))^{m\times m}$, $\mathbf{f}\in C((0,T); \mathcal{H}^s(M))^m$, $\mathbf{v}_0\in\mathcal{H}^s(M)^m$, where $s>\frac{N}{2}+1$ is an integer. Then there exists a unique solution to \eqref{ar44}, \eqref{ar45}, i.e. a function 
\begin{displaymath}
\mathbf{v}\in \left[C([0,T); \mathcal{H}^s(M))\cap C^1((0,T); \mathcal{H}^{s-1}(M))\right]^m
\end{displaymath}
satisfying \eqref{ar44} and \eqref{ar45} pointwise (i.e. in the classical sense).
\end{theorem}

\begin{definition}\label{equic2}
Let $X$ be a compact metric space, $Y$ a Banach space, and $C(X,Y)$ the Banach space of continuous functions from $X$ to $Y$ with the sup norm. A subset $\Omega$ of $C(X,Y)$ is equicontinuous if, for every
$\epsilon > 0$, there is a $\delta > 0$ such that, for each $f\in\Omega$, $d(x,y) < \delta$ implies
\begin{equation}
\| f(x)-f(y)\|\leq\epsilon,
\end{equation}
where $d$ is the metric in $X$. \\
$\Omega$ is uniformly bounded or equibounded if there exists $M\in\mathbb{R}^{+}$ such that 
\begin{equation}
\|f(x)\|\leq M \ \ \ \ \ \ \ \ \left(x\in X; f\in\Omega\right).
\end{equation}
\end{definition}

\subsection{Statement of the main results}

In this section we state the main results of this paper. The local existence and uniqueness of a solution for Models A1 and A2 is stated in the following theorems.

\begin{theorem}[\textbf{Model A1}]\label{TeoremaBr1}
Let us consider the system
\begin{eqnarray}\label{Tbr1}
\begin{cases}
\partial_t\eta(t,z)=\zeta(t,z),\\
\widehat{\alpha}A(t,z)\partial_{t}\zeta(t,z)+\widetilde{k}\zeta(t,z)+\kappa\eta(t,z)-A(t,z)P(t,z)+A(t,z)\widetilde{P}=0, \\
\left(\widetilde{h}+a(t)+\eta(t,z)\right)\partial_t A(t,z)+\left(a'(t)+\zeta(t,z)+u(t,z)\right)A(t,z)-Q_p=0,\\
\rho\partial_t u(t,z)+\rho u(t,z)\partial_z u(t,z)+\beta u(t,z)+\partial_z P(t,z)=0,  \\
R\partial_t P(t,z)-KP(t,z)^2-KP(t,z)\left(RQ_p+\widetilde{P}\right)=0,
\end{cases}
\end{eqnarray}
where $t\in[0,T_0]$, $z\in[0,L]$, $a(t)\in C^{\infty}([0,T_0])$, with the initial data  
\begin{align}\label{Tbr2}
u(0,z)=u_0(z), &\ \ \eta(0,z)=\eta_0(z), \ \ \zeta(0,z)=\zeta_0(z), \notag\\
 P(0,z)&=P_0(z), \ \ A(0,z)=A_0(z),
\end{align}
in $\mathcal{H}^s([0,L])$ with $\displaystyle{s>\frac{5}{2}}$. Consider the boundary conditions \eqref{br4}, \eqref{br11} and \eqref{br15} for the system \eqref{Tbr1} and assume that the compatibility conditions \eqref{br6}, \eqref{br12} and \eqref{br16} are satisfied.\\
Then, there exists a time $T<T_0$ such that the problem \eqref{Tbr1} with the previous initial and boundary conditions has a local unique solution
\begin{displaymath}
\mathcal{X}(t,z)=\left(u(t,z), \eta(t,z), \zeta(t,z), P(t,z), A(t,z)\right),
\end{displaymath}
such that
\begin{displaymath}
\mathcal{X}(t,z)\in \left[C^1((0,T)\times [0,L])\right]^5.
\end{displaymath}
\end{theorem}

\begin{theorem}[\textbf{Model A2}]\label{TeoremaBr2}
Let us consider the system
\begin{eqnarray}\label{Tbr21}
\begin{cases}
\partial_t\eta(t,z)=\zeta(t,z),\\
\widehat{\alpha}A(t,z)\partial_{t}\zeta(t,z)+\widetilde{k}\zeta(t,z)+\kappa\eta(t,z)-A(t,z)P(t,z)+A(t,z)\widetilde{P}=0, \\
\left(\widetilde{h}+a(t)+\eta(t,z)\right)\partial_t A(t,z)+\left(a'(t)+\zeta(t,z)\right)A(t,z)-Q_p\\
\ \ \ \ \ \ \ \ \ \ \ \ \ \ \ \ \ \ \ \ \ \ \ \ \ \ \ \ \ \ \ \ \ \ \ \ \ \ \ \ \ \ \ \ \ \ \ \ \ \ \ \ \ \ \ \ \ \ \ \ \  +\displaystyle{\frac{1}{R}}\left(P(t,z)-\widetilde{P}\right)=0, \\
\rho\partial_t u(t,z)+\rho u(t,z)\partial_z u(t,z)+\beta u(t,z)+\partial_z P(t,z)=0,  \\
R\partial_t P(t,z)-KP(t,z)^2-KP(t,z)\left(RQ_p+\widetilde{P}\right)=0,
\end{cases}
\end{eqnarray}
where $t\in[0,T_0]$, $z\in[0,L]$, $a(t)\in C^{\infty}([0,T_0])$, with initial conditions \eqref{Tbr2}. Consider the boundary conditions \eqref{br4}, \eqref{br11} and \eqref{br15} for the system \eqref{Tbr1} and assume that the compatibility conditions \eqref{br6}, \eqref{br12} and \eqref{br16} are satisfied.\\
Then, there exists a time $T<T_0$ such that the problem \eqref{Tbr21} has a local unique solution
\begin{displaymath}
\mathcal{X}(t,z)=\left(u(t,z), \eta(t,z), \zeta(t,z), P(t,z), A(t,z)\right),
\end{displaymath}
such that
\begin{displaymath}
\mathcal{X}(t,z)\in \left[C^1((0,T)\times [0,L])\right]^5.
\end{displaymath}
\end{theorem}
By following the same lines of arguments of the Theorem 2 stated in \cite{donmaro}, we are able to prove global existence and uniqueness of solutions for Models A1 and A2 under some restriction conditions on the initial data. \\
For the sake of simplicity, henceforth we will specify in brackets the system and the conditions related to Model A2 when it admits the same results of Model A1.\\
The global existence and uniqueness theorem can be stated as follows 
\begin{theorem}[\textbf{Model A1 (Model A2)}]\label{TeoremaBr3}
Let us consider the system \eqref{Tbr1} (\eqref{Tbr21}) with initial conditions \eqref{Tbr2} and boundary conditions \eqref{br4}, \eqref{br11} and \eqref{br15}. Assume that the compatibility conditions \eqref{br6}, \eqref{br12} and \eqref{br16} are satisfied and that
\begin{equation}\label{Tbr5}
f'(z)\geq-\frac{\beta}{\rho}, 
\end{equation}
where $f(z)=u(0,z)$, and 
\begin{equation}\label{Tbr5bis}
\|b(z)\|_{\mathcal{H}^{s}_z}\leq\frac{\beta}{4\widehat{\mathcal{C}}_1\rho},
\end{equation}
 if $b(z)<0$ with $b(z)=P(0,z)$, and
\begin{equation}\label{Tbr5tris}
\|b(z)\|_{\mathcal{H}^{s}_z}\leq\varepsilon\frac{\beta}{4\widehat{\mathcal{C}}_1\rho},
\end{equation}
if $b(z)>0$, where $\widehat{\mathcal{C}}_1=\widehat{\mathcal{C}}_1\left(\|f(z)\|_{\mathcal{H}^s_z}\right)$ and $\varepsilon=\varepsilon(\widetilde{T}_0)$ are constants.\\
Then there exists a global unique solution 
\begin{displaymath}
\mathcal{X}(t,z)=\left(u(t,z), \eta(t,z), \zeta(t,z), P(t,z), A(t,z)\right),
\end{displaymath}
to the problem \eqref{Tbr1},\eqref{Tbr2} (\eqref{Tbr21},\eqref{Tbr2}) such that
\begin{displaymath}
\mathcal{X}(t,z)\in \left[C^1((0,T]\times[0,L] )\right]^5,
\end{displaymath}
for every $T\geq0$.
\end{theorem}
The proofs of Theorems \ref{TeoremaBr1}, \ref{TeoremaBr2} and \ref{TeoremaBr3} are treated in detail in the Appendix \ref{appendix}.

\section{NUMERICAL SIMULATIONS}

In the present section, direct numerical simulations are carried out to investigate the accuracy and completeness of the theorems we proved. \\
In this study, the CSF is treated as a Newtonian fluid with a dynamic viscosity and density equal to $1003\cdot10^{-3}$ $\kg/(\meter\cdot\second)$ and $998,2$ $\kg/\meter^3$, respectively \cite{linninger2007}. The brain tissue is considered as a linear viscoelastic material with the storage and loss moduli of $2038$ and $1356\,\Pa$ for the healthy subjects, and the density of $1040$ $\kg/\meter^3$ (\cite{ghol}, \cite{streit}). The CSF flow rate in the lateral ventricles is $0.35$ $\cm^3/\min$ \cite{pople}. This value is used as the amplitude in the input fluid pulsatile flow rate function for numerical models. The final section of the ventricular system after the fourth ventricle is selected for the flow output location in Model A2. We set the normal baseline CSF pressure at $500$ $\Pa$. \\
The computation is performed in 2D under the assumption of axial symmetry of the problems. Models A1 and A2 are discretized by taking into account the approximating systems defined in Section \ref{approx} and are solved by using the Runge-Kutta MATLAB solver ODE45, with the initial mesh size determined by $dz=L/nz$, where $L$ is fixed equal to 1 for the sake of simplicity, and $nz$ is the total number of meshes that we take equal to 100. Finally we set the time discretization with $\Delta t = 5 \times 10^{-3}$.\\ 
A detailed analysis of the numerical methods we shall adopt for the models analyzed is beyond the scope of the present paper, and will be subject of future investigation.\\

The present section is organized as follows.
\begin{itemize}
\item First of all, we prescribe the boundary conditions \eqref{br4}, \eqref{br11} and \eqref{br15} for $u, \eta, \zeta,$ while by taking into account \eqref{br18} and \eqref{br26} the boundary conditions for the pressure and the axial section are the following 
\begin{align}\label{br21}
P(t,0)&=\displaystyle{\frac{b(0)e^{\mathcal{C}t}}{1+\displaystyle{\frac{K}{R\mathcal{C}}}b(0)\left(1-e^{\mathcal{C}t}\right)}},\notag \\
P(t,L)&=\displaystyle{\frac{b(L)e^{\mathcal{C}t}}{1+\displaystyle{\frac{K}{R\mathcal{C}}}b(L)\left(1-e^{\mathcal{C}t}\right)}}, \notag \\
A(t,0)&=h(0)e^{G(t,0)t}+e^{-G(t,0)t}\int_0^t H(s,0)e^{G(s,0)s}\,ds,\notag \\
A(t,L)&=h(L)e^{G(t,L)t}+e^{-G(t,L)t}\int_0^t H(s,L)e^{G(s,L)s}\,ds, 
\end{align}
where 
\begin{equation}\label{br30}
G(t,0)=\displaystyle{\frac{\displaystyle{a'(t)}}{\widetilde{h}+a(t)}}, \ \ \ \ H(t,0)=\displaystyle{\frac{Q_p}{\widetilde{h}+a(t)}},
\end{equation}
in Model A1 and
\begin{equation}\label{br30Bis}
G(t,0)=\displaystyle{\frac{\displaystyle{a'(t)}}{\widetilde{h}+a(t)}}, \ \ \ \ H(t,0)=\displaystyle{\frac{RQ_p-P(t,0)-\widetilde{P}}{R(\widetilde{h}+a(t))}},
\end{equation}
in Model A2. \\

Then we run numerical simulations in order to validate the main Theorem of the paper which proves global existence of solutions provided that proper conditions on the initial data of the CSF flow velocity and of the intracranial pressure are satisfied. We perform this task by assuming 
\begin{align}\label{numsim1}
u_0(z)&=4\sin(\pi z)+1, \notag \\
P_0(z)&=\displaystyle{\frac{1}{6}}\cos(\pi z),
\end{align}
which are not so far from a good real approximation of the CSF flow velocity in a small compartment perfused by the fluid. \\
As a byproduct of our numerical results, we obtain interesting informations about the behaviours of the quantities involved in these models.
\item As a second task, we implement again the models numerically but we fix initial data $u_0(z)$ and $P_0(z)$ which violate the conditions \eqref{Tbr5} and \eqref{Tbr5bis} (we show here only the case of initial intracranial pressure in upright posture) of the Theorem \ref{TeoremaBr3}. We shall compare these simulations to the analysis developed in this paper and show that our approach is reliable. In order to do that we assume
\begin{align}\label{numsim2}
u_0(z)&=-(\exp(z)+1), \notag \\
P_0(z)&=\exp(z).
\end{align}
\end{itemize}
\vspace{2pt}
\begin{flushleft}
\textbf{First case:} \textit{the initial data fulfill the conditions of the Theorem \ref{TeoremaBr3}.}\label{caseA}
\end{flushleft}
\begin{figure}[htbp]
\centering
\includegraphics[width=\textwidth]{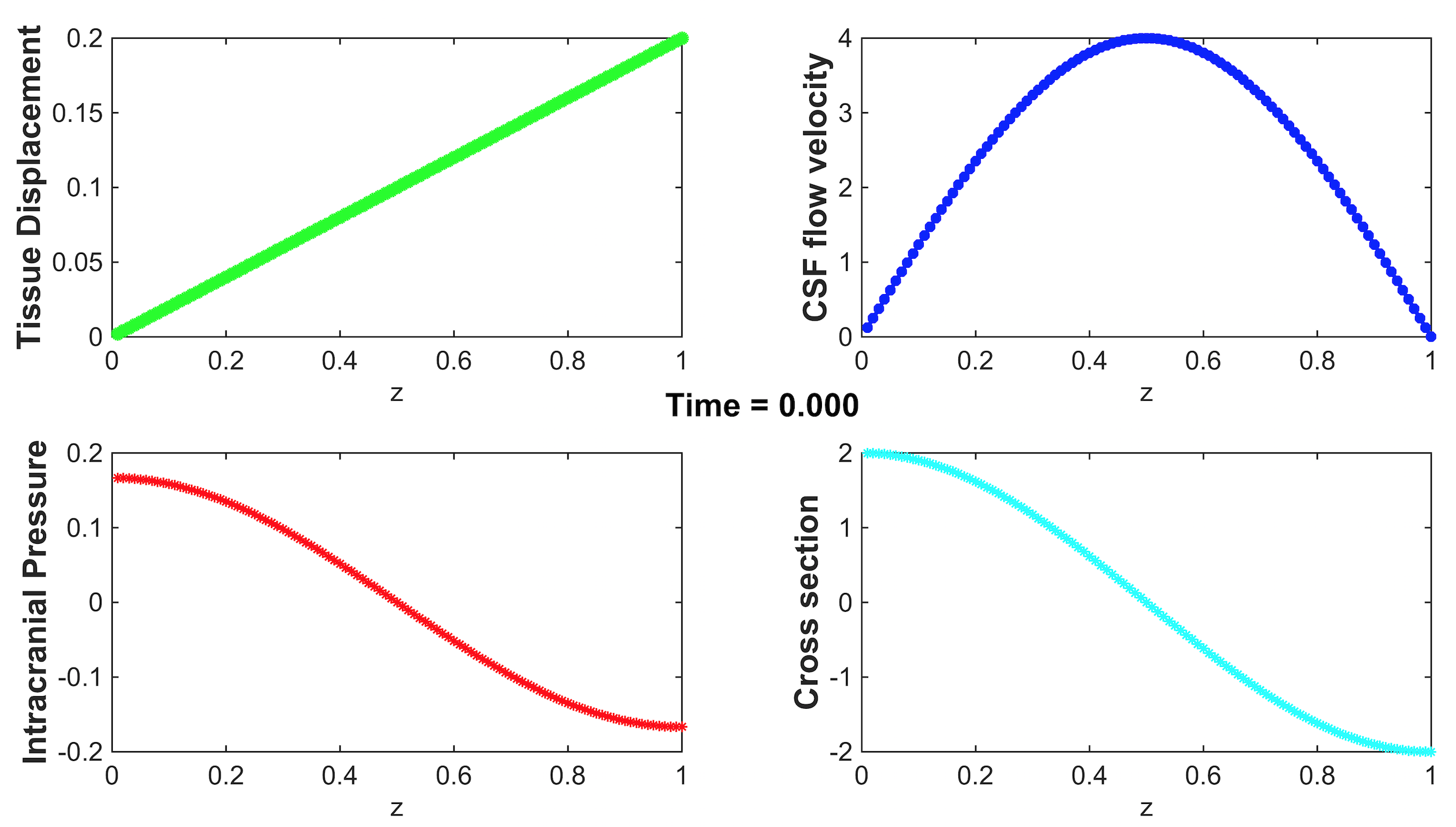}
\caption{The initial configuration of $P, \eta, u, \zeta $ and $A$ for Model A1.\label{FiguraSim1}}
\end{figure}
\begin{figure}[htbp]
\centering
\includegraphics[width=\textwidth]{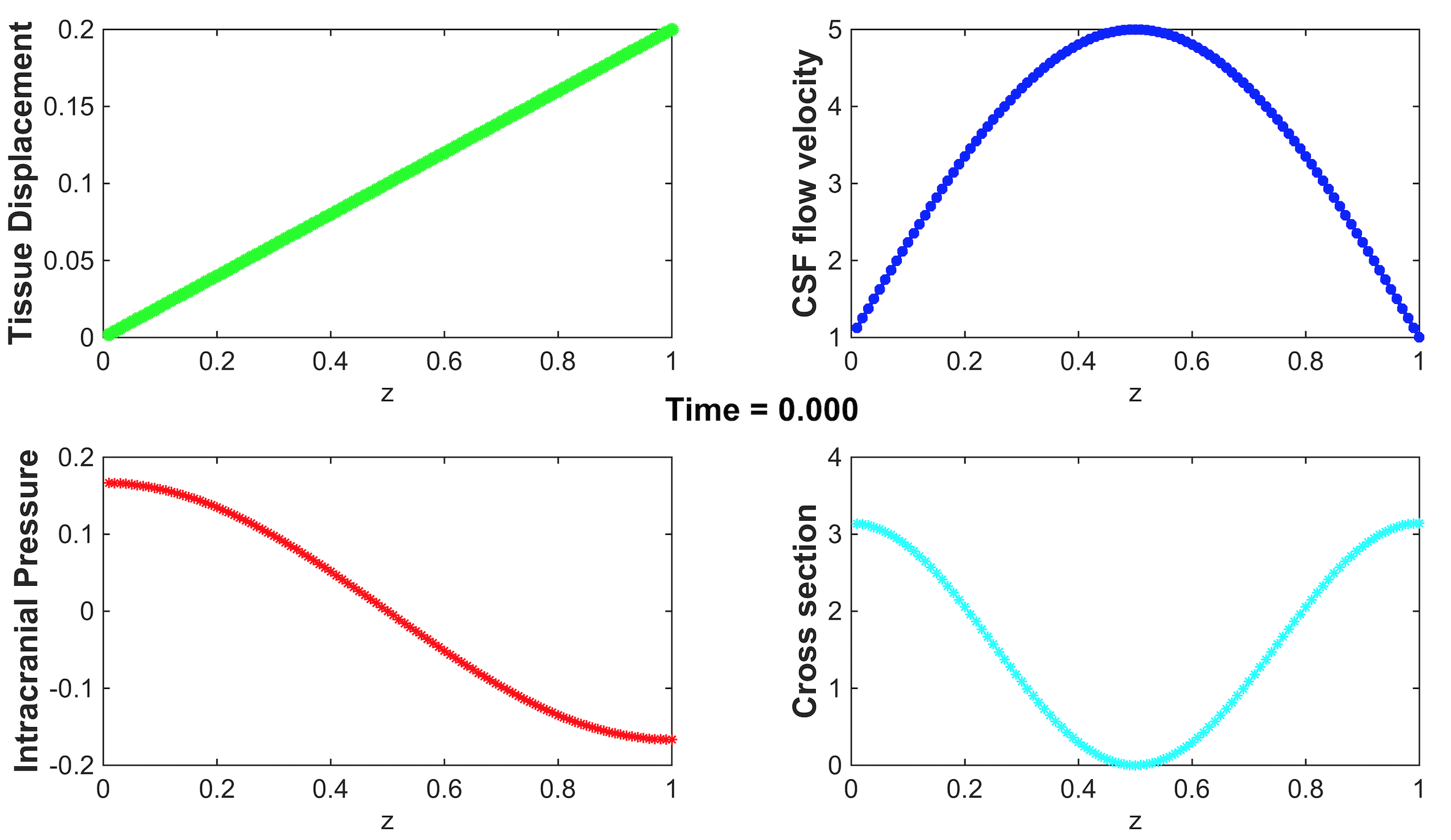}
\caption{The initial configuration of $P, \eta, u, \zeta $ and $A$ for Model A2.\label{FiguraSim2}}
\end{figure}
We perform the simulation by assuming the initial condition \eqref{numsim1} and arbitrary data for the others quantities for which we choose
\begin{align}\label{numsim3}
\eta_0(z)&=\frac{1}{5}z,\notag \\
 \zeta_0(z)&=\frac{1}{2}z+1,\notag \\
 A_0(z)&=2\cos(\pi z).
\end{align}
The initial configurations of Model A1 and A2 are displayed in Figures \ref{FiguraSim1} and \ref{FiguraSim2}.\\
\begin{figure}[htbp]
\centering
\includegraphics[width=\textwidth]{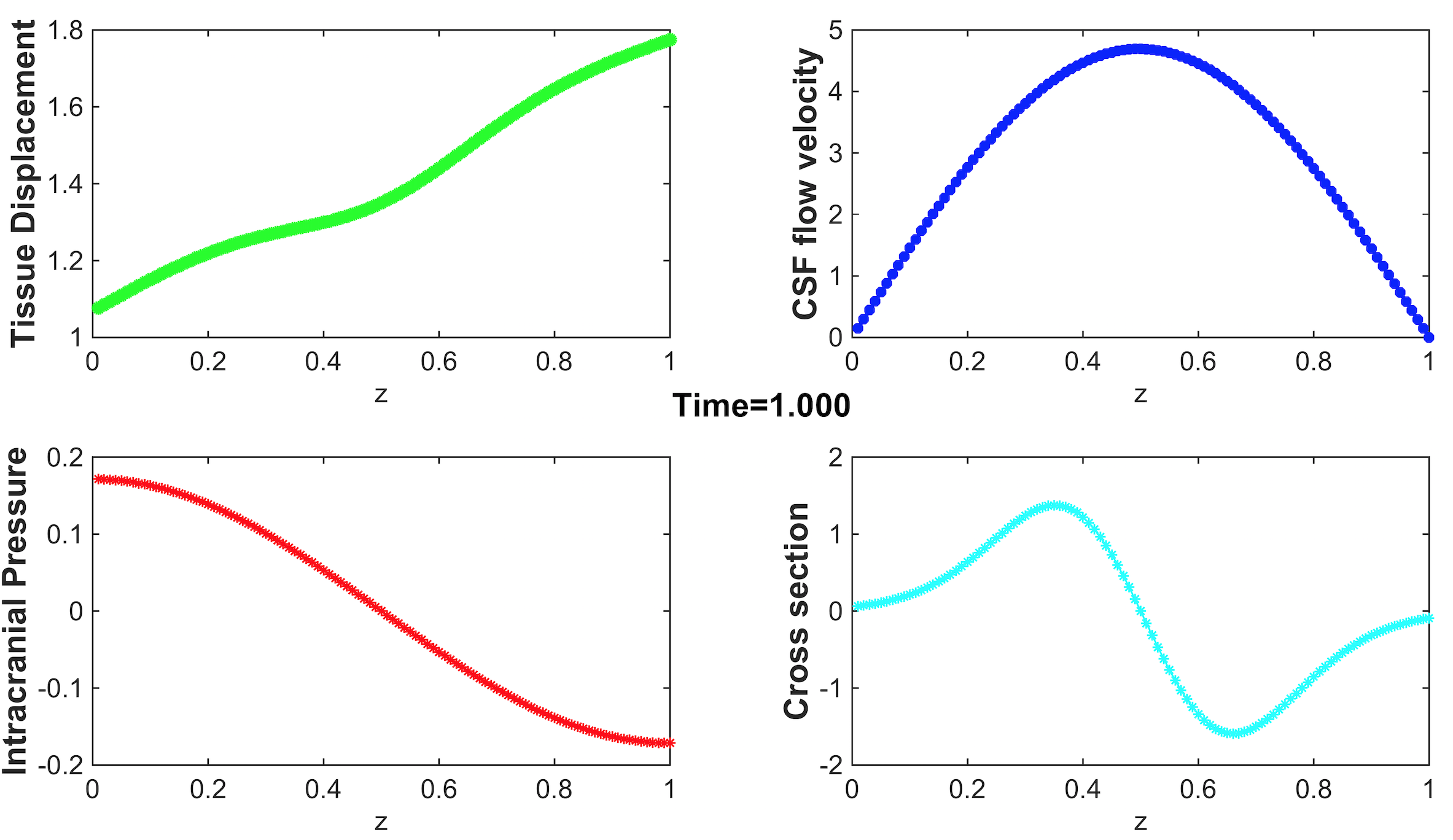}
\caption{Evolution of $P, \eta, u, \zeta $ and $A$ for Model A1.\label{FiguraSim3}}
\end{figure}
\begin{figure}[htbp]
\centering
\includegraphics[width=\textwidth]{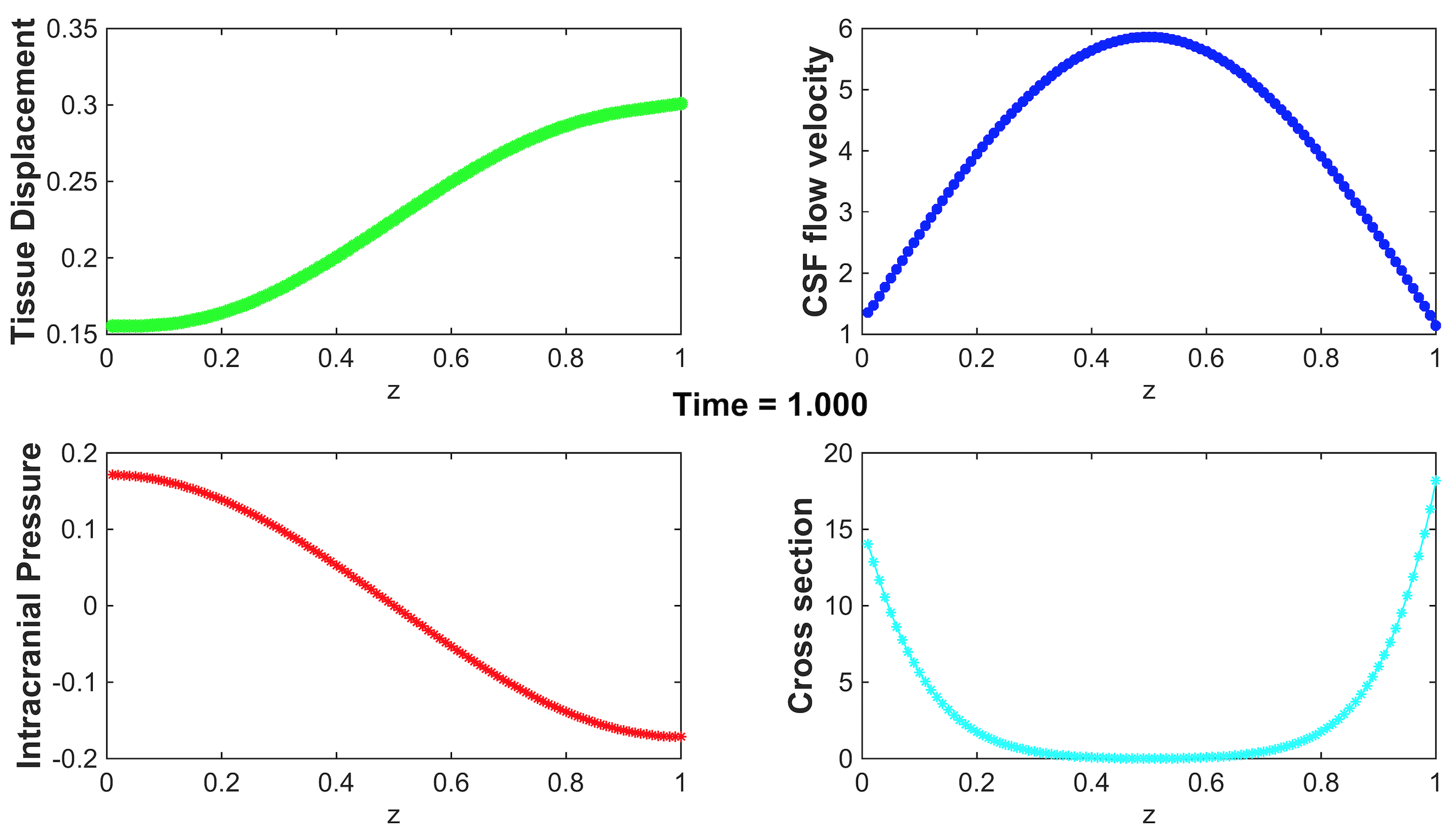}
\caption{Evolution of $P, \eta, u, \zeta $ and $A$ for Model A2.\label{Sim4}}
\end{figure}
In Figures \ref{FiguraSim3} and \ref{Sim4} we can observe the evolution of all the quantities involved in the CSF dynamics in a single compartment of length $L$ and after a certain time $\bar{t}>0$.\\
\begin{itemize}
\item In \textbf{Model A1}:
\begin{itemize}
\item the intracranial pressure, $P(t,z)$, shows the maximum value at the inlet and then it decreases through the compartment;
\item the tissue displacement, $\eta(t,z)$ undergoes a small decrease close to the midsection and an increase at the boundary, which correspond to a compression and an expansion of the compartment respectively;
\item the flow velocity, $u(t,z)$, shows a parabolic profile that represents exactly what we expected since the CSF flow is laminar; moreover, by taking into account the initial configuration, it is possible to notice that at time $\bar{t}$ the velocity at midsection is decreased;
\item for the cross section, $A(t,z)$, we can observe an increase, then an enlargement of the section, in the first part of the compartment and a narrowing in the second part.
\end{itemize}
Therefore, corresponding to the maximum value of the pressure we have an increase in cross section which shrinks in relation to the decrease of the intracranial pressure. In fact, it occurs a phenomenon called Venturi effect for which in a tube with variable section the pressure decreases as section decreases (by the continuity equation). Strangely, the compression and expansion of the tissue in the compartment does not reflect the behavior of the cross section, in fact a tissue compression occurs in relation to the increase of the section and vice versa, while we would expect the opposite. 
\item In \textbf{Model A2}
\begin{itemize}
\item the intracranial pressure shows the same behavior as in Model A1; this is due to the fact that the Riccati equation which really rules the pressure evolution in time, in general seems to evolves independently of the other quantities involved in the models;
\item the tissue displacement increases continuously through the compartment by starting from a minimum value at the inlet;
\item the flow velocity behaves as in Model A1 but the value at midsection in this model is greater than the corresponding value in the initial configuration;
\item cross section decreases up to a minimum value at the midsection and in the second part of the compartment it shows an increase.
\end{itemize}
Thus, as in Model A1, in Model A2 an unusual behavior occurs when the intracranial pressure reaches its minimum value at the outlet since the cross section shows an increase rather than a compression. On the other hand, tissue displacement and cross section behaviors are properly related, in particular at the outlet of the compartment, where the maximum value of $\eta(t,z)$ corresponds to a section enlargement.
\end{itemize} 
The behaviour of the parameters analyzed here, up to the final time $\bar{t}=1.000$, is due to the boundary conditions we imposed for both models and to the fact we are neglecting fundamentals interactions in the cerebral pattern, that are assumptions which do not reflect the real CSF physiology and represent a limitation in the intracranial detailed description. \\
But the main result obtained with these first simulations is the regularity of the solutions displayed above which is in good agreement with the Theorem \ref{TeoremaBr3}.

\vspace{12pt}
\begin{figure}[h]
\centering
\includegraphics[width=\textwidth]{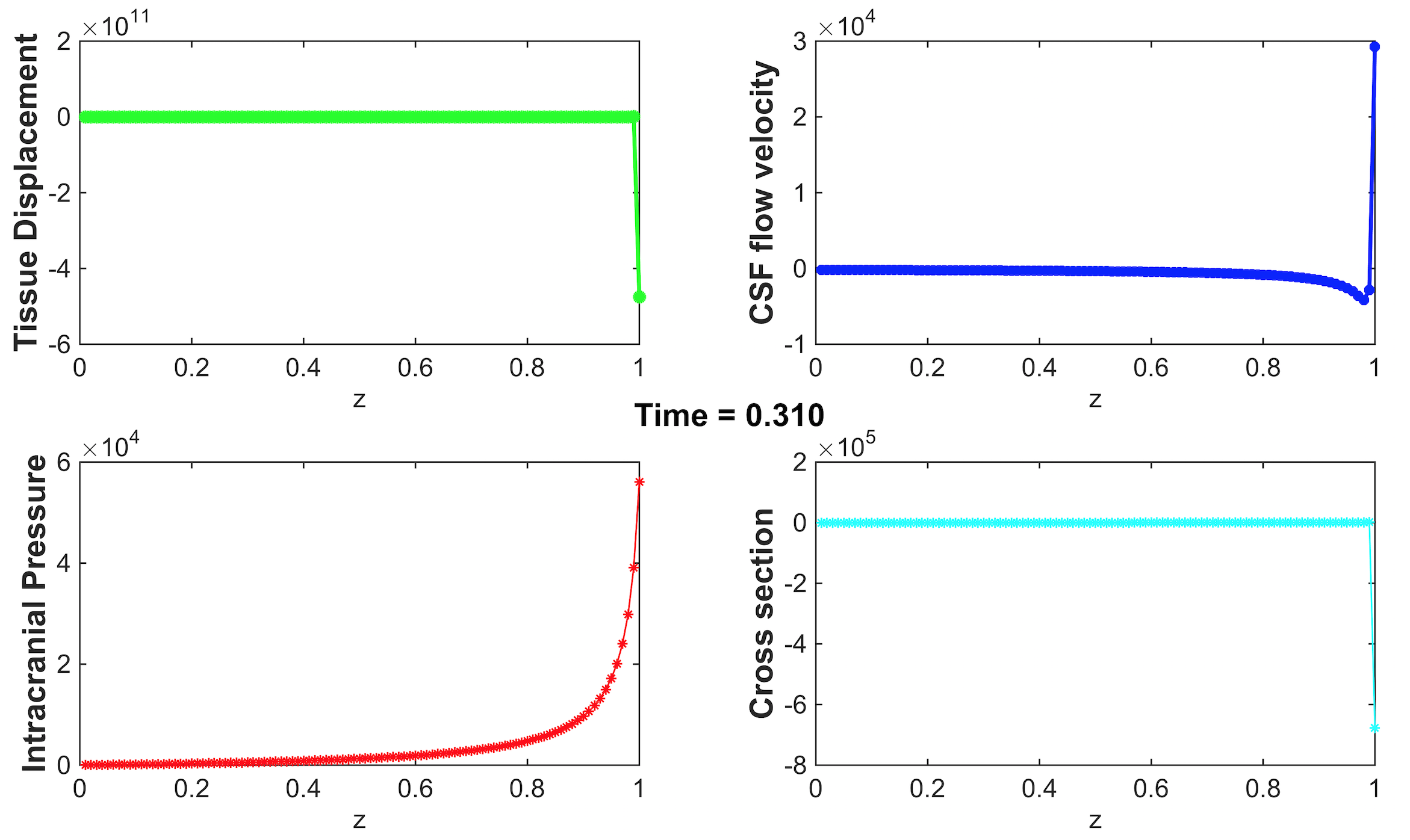}
\caption{Blow-up which occurs in Model A1 with initial conditions \eqref{numsim2}, \eqref{numsim3}.\label{BU3}}
\end{figure}
\newpage
\begin{flushleft}
\textbf{Second case:} \textit{the initial data violate the conditions of the Theorem \ref{TeoremaBr3}.}
\end{flushleft}

\begin{figure}[h]
\centering
\includegraphics[width=\textwidth]{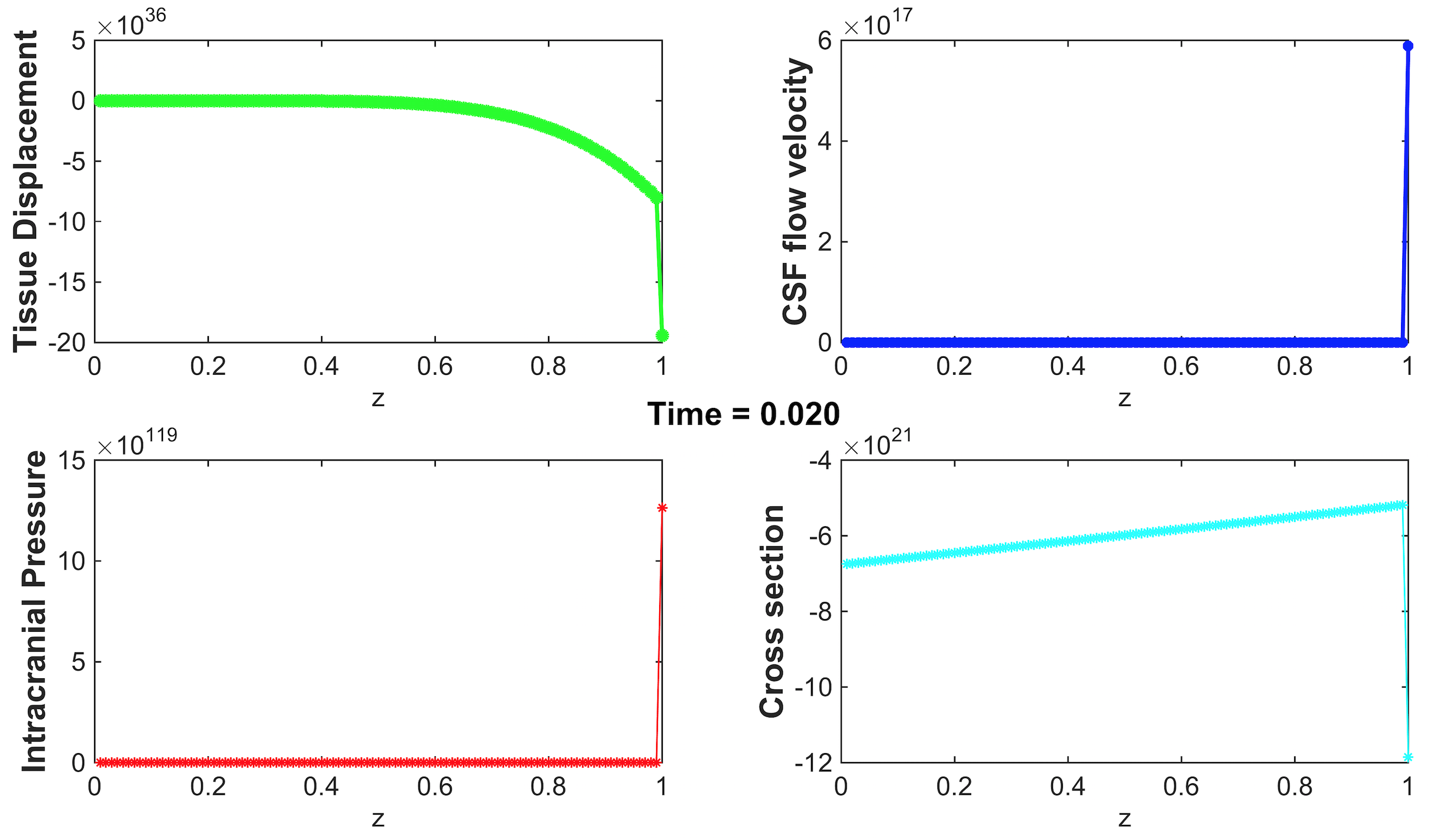}
\caption{Blow-up which occurs in Model A2 with initial conditions \eqref{numsim2}, \eqref{numsim3}.\label{BU4}}
\end{figure}

Figures \ref{BU3} and \ref{BU4} show precisely what we predicted with our analysis: when the global existence conditions \eqref{Tbr5} and \eqref{Tbr5bis} (or \eqref{Tbr5tris}) are violated, the CSF velocity flow and the pressure blow up in a finite time and the other quantities, which inherit the same behaviour, exhibit simultaneously a blow up.\\
The occurrence of blow up underlines the fact that we are neglecting important elements for the CSF dynamics, in fact in nature it is impossible to detect this kind of effect in a human body. 

\section{Comparison to real data}\label{comparison}
Once validated the main results stated in the present paper, a comparison with respect to the real dynamics detected in the intracranial pattern is required in order to understand how mathematical assumptions work and are able to alter the effective behavior of all quantities involved in both Models A1 and A2.\\ 
\vspace{6pt}
\begin{table}[h]
\centering
\begin{tabular}{|r|c|c|}
\hline
\textbf{Clinical observation}&\textbf{Value}&\textbf{Reference}\\ \hline
CSF flow velocity& $50-80\ \ \mm/\s$ &\cite{sharma}\\ \hline
ICP (supine position) & $7 -15\ \ \mmHg$  &\cite{alperin}\\ \hline
ICP (upright posture) & average value of $-3,4\ \ \mmHg$  &\cite{alperin} \\
\hline
\end{tabular}
\caption{Published physiological data.\label{tabellina}}
\end{table}
\begin{figure}[htbp]
\centering
\includegraphics[width=\textwidth]{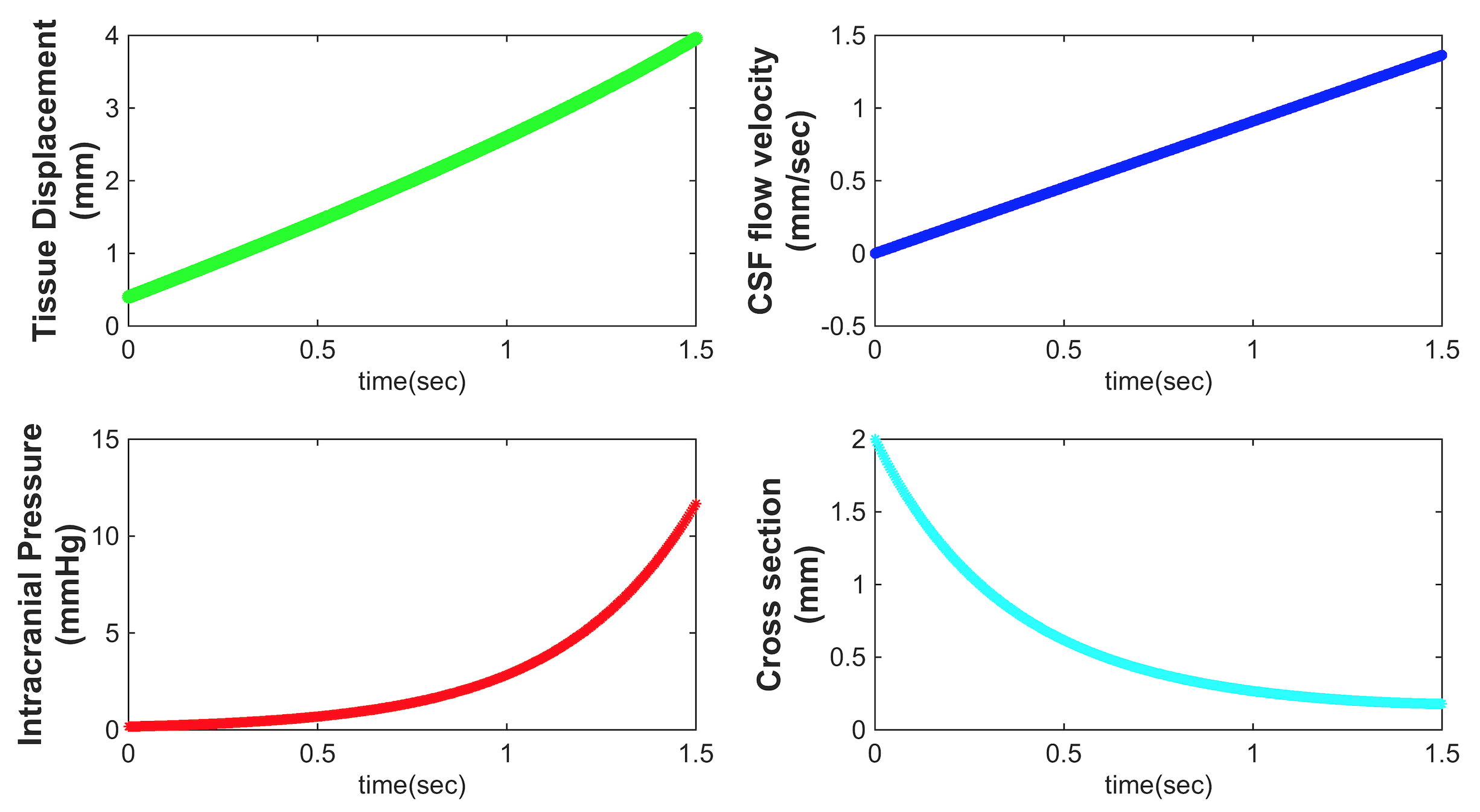}
\caption{Evolution in time of Model A1.\label{EvolTemp1}}
\end{figure}
\begin{figure}[htbp]
\centering
\includegraphics[width=\textwidth]{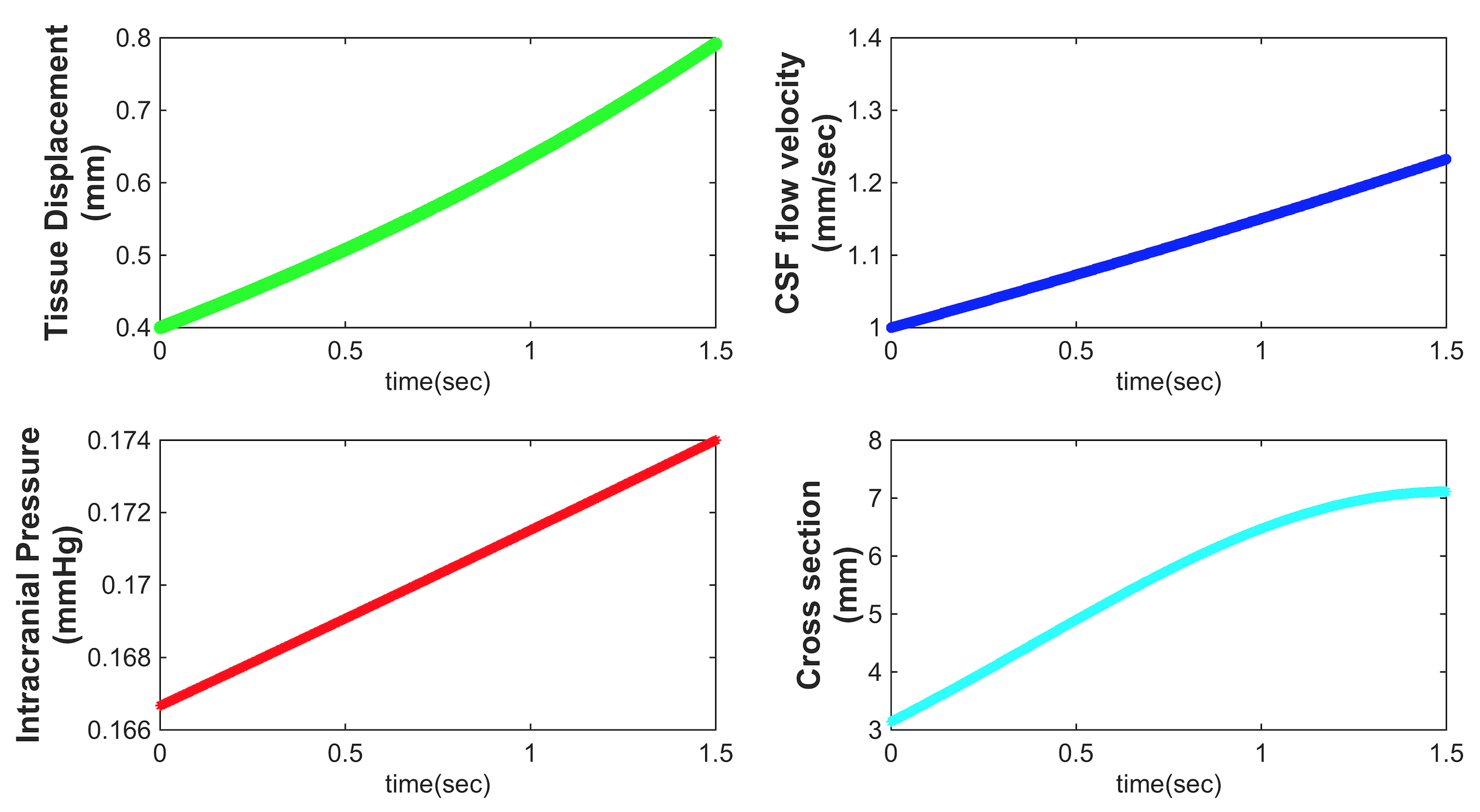}
\caption{Evolution in time of Model A2.\label{EvolTemp2}}
\end{figure}

It is interesting to compare in particular the CSF flow velocity and the ICP values available in literature, thanks to the modern MRI techniques, with the results obtained in our analysis and simulations. In Fig.\,\ref{EvolTemp1} we can observe that in Model A1, with a peak of $1,4 \ \ \mm/\s$, CSF flow velocity is very slow with respect to the value indicated in Table \ref{tabellina}. On the contrary, even though the boundary conditions adopted here, the values assumed by the intracranial pressure are quite consistent with the real data since after the first second it lies in the range of values specified in Table \ref{tabellina}. \\

In Fig.\,\ref{EvolTemp2} the plot of the CSF flow velocity shows a behavior similar to the one displayed in Fig.\,\ref{EvolTemp1} but the maximum value is about $1,23\ \ \mm/\s$ which is smaller than the value detected in Model A1, as well as the value listed in Table \ref{tabellina}. This is due to the fact that Model A2 describes the reabsorption process of CSF characterized by a velocity gradient decrease and by a damping effort of the resistance to the CSF absorption. The intracranial pressure is affected by a negligible alteration and the evolution plotted is very far from the data expected in the real physiology either in a supine position or in the upright posture.\\

Concerning the tissue displacement and the cross section, the comparison with real data is not self-evident since the compartments perfused by the CSF are very different each other either in anatomy or in physiology. Moreover, Models A1 and A2 do not take into account important interactions among the systems involved in the CSF dynamics, therefore a more detailed CSF model (see \cite{linninger2007}) would be required in order to have a very accurate discussion regarding the differences pointed out between model predictions and published data.

\section{Conclusions}
This paper presents an extensive analysis of two CSF models, Models A1 and A2, that represent an improvement with respect to the study provided in \cite{donmaro}.\\
As a first attempt in this direction, starting by the CSF model proposed by Linninger et al.\,in 2005 (see \cite{b3}) we have adopted a simplified model that describes the CSF flow without considering the different and complex interactions that occur inside the skull. We have also included to both models an intracranial pressure equation with the aim of ensuring a more detailed representation than the one showed in \cite{donmaro}, although it provides a cruder representation with respect to more recent CSF models.\\
This choice is intended to be only the first step towards the analysis of much more realistic and anatomically complete models (see \cite{linninger2009}), able to describe every interaction that occurs between the three main subsystems in the intracranial pattern: cerebral vasculature, CSF, porous parenchyma.\\

In the present work we first demonstrated the local existence and uniqueness in time of solutions to both Models A1 and A2, then we proved the global existence and uniqueness of solutions under stringent conditions on the initial data.\\
Furthermore, the numerical simulations supported the theoretical results and show evidence that the proposed analysis is consistent as a starting point also in cases in which more parameters and physiological quantities are involved in the models in order to guarantee a more comprehensive approach to the CSF dynamics. The authors will continue the study of this problem in a future work.

\appendix

\section*{APPENDIX}\label{appendix}

This appendix includes the detailed proofs of the main theorems stated in the present paper.
\vspace{12pt}
\section{PROOF OF THE THEOREMS \ref{TeoremaBr1} AND \ref{TeoremaBr2}}

In this Section we are going to prove the first result of our paper. In order to do that we need to proceed in two steps: we set up an iterative process by means of approximating systems associated to \eqref{br1bis} and \eqref{br2bis} and then we prove the convergence of the approximating solutions to classical solutions of Models A1 and A2.

\subsection{Approximating scheme}\label{approx}

For the purpose of being more precise we define one by one the approximations systems of \eqref{br1bis} and \eqref{br2bis} as follows.

\begin{flushleft}
\textbf{Approximating system A1 (As A1)}
\end{flushleft}
\begin{eqnarray}\label{br1bbis}
\begin{cases}
\partial_t\eta^{n+1}(t,z)=\zeta^{n}(t,z),\\
\widehat{\alpha}A^{n}(t,z)\partial_{t}\zeta^{n+1}(t,z)\!+\!\widetilde{k}\zeta^{n+1}(t,z)\!+\!\kappa\eta^{n}(t,z)\!-\!A^{n}(t,z)P^{n}(t,z)\!+\!A^{n}(t,z)\widetilde{P}\!=\!0, \\
\left(\widetilde{h}\!+\!a(t)\!+\!\eta^{n}(t,z)\right)\partial_t A^{n+1}(t,z)\!+\!\left(a'(t)\!+\!\zeta^{n}(t,z)\!+\!u^{n}(t,z)\right)A^{n+1}(t,z)\!-\!Q_p\!=\!0,\\
\rho\partial_t u^{n+1}(t,z)+\rho u^n(t,z)\partial_z u^{n+1}(t,z)+\beta u^{n+1}(t,z)+\partial_z P^n(t,z)=0,  \\
R\partial_t P^{n+1}(t,z)-KP^{n+1}(t,z)^2-KP^{n+1}(t,z)\left(RQ_p+\widetilde{P}\right)=0,
\end{cases}
\end{eqnarray}
\newpage
\begin{flushleft}
\textbf{Approximating system A2 (As A2)}
\end{flushleft}
\begin{eqnarray}\label{br2bbis}
\begin{cases}
\partial_t\eta^{n+1}(t,z)=\zeta^n(t,z),\\
\widehat{\alpha}A^{n}(t,z)\partial_{t}\zeta^{n+1}(t,z)\!+\!\widetilde{k}\zeta^{n+1}(t,z)\!+\!\kappa\eta^n(t,z)\!-\!A^n(t,z)P^n(t,z)\!+\!A^{n}(t,z)\widetilde{P}\!=\!0, \\
\left(\widetilde{h}+a(t)+\eta^n(t,z)\right)\partial_t A^{n+1}(t,z)+\left(a'(t)+\zeta^n(t,z)\right)A^{n+1}(t,z)-Q_p\\
\ \ \ \ \ \ \ \ \ \ \ \ \ \ \ \ \ \ \ \ \ \ \ \ \ \ \ \ \ \ \ \ \ \ \ \ \ \ \ \ \ \ \ \ \ \ \ \ \ \ \ \ \ \ \ \ \ \ \ \ \ \ \ \ \ +\displaystyle{\frac{1}{R}}\left(P^n(t,z)-\widetilde{P}\right)=0, \\
\rho\partial_t u^{n+1}(t,z)+\rho u^n(t,z)\partial_z u^{n+1}(t,z)+\beta u^{n+1}(t,z)+\partial_z P^n(t,z)=0,  \\
R\partial_t P^{n+1}(t,z)-KP^{n+1}(t,z)^2-KP^{n+1}(t,z)\left(RQ_p+\widetilde{P}\right)=0,
\end{cases}
\end{eqnarray}
where $z\in[0,L]$, $t\in(0, T_0)$, with $T_0\geq T>0$. We assume for both systems that the initial conditions
\begin{align}\label{br3bbis}
u(0,z)^{n+1}=f(z), \ \ & \eta(0,z)^{n+1}=g(z), \ \ \zeta(0,z)^{n+1}=\zeta_0(z)=q(z),\notag \\
P(0&,z)^{n+1}=b(z),\ \ A(0,z)^{n+1}=h(z),
\end{align}
while the boundary conditions for both systems are as follows:
\begin{align}\label{br4bbis}
u(t,0)^{n+1}&=u(t,L)^{n+1}=0, \ \ \eta(t,0)^{n+1}=\eta(t,L)^{n+1}=0, \notag\\
&\zeta(t,0)^{n+1}=\zeta(t,L)^{n+1}=0.
\end{align}
For both systems \eqref{br1bbis} and \eqref{br2bbis} we fix the initial step as,
\begin{align}\label{br8bbis}
u(0,z)^0=u_0(z)&=f(z), \ \ \eta(0,z)^0=\eta_0(z)=g(z), \ \ \zeta(0,z)^0=\zeta_0(z)=q(z),\notag \\
P(0,z)^0&=P_0(z)=b(z),\ \ A(0,z)^0=A_0(z)=h(z),
\end{align}
and we proceed by induction on the iteration number $n$.

\begin{flushleft}
\textbf{BASIC STEP: $\mathbf{n=0}$.}
\end{flushleft}
The approximating systems take the following form
\begin{flushleft}
\textbf{As A1}
\end{flushleft}
\begin{eqnarray}\label{br9bbis}
\begin{cases}
\partial_t\eta^{1}(t,z)=\zeta^{0}(t,z),\\
\widehat{\alpha}A^{0}(t,z)\partial_{t}\zeta^{1}(t,z)+\widetilde{k}\zeta^{1}(t,z)+\kappa\eta^{0}(t,z)-A^{0}(t,z)P^{0}(t,z)+A^{0}(t,z)\widetilde{P}=0, \\
\left(\widetilde{h}\!+\!a(t)\!+\!\eta^{0}(t,z)\right)\partial_t A^{1}(t,z)\!+\!\left(a'(t)\!+\!\zeta^{0}(t,z)\!+\!u^{0}(t,z)\right)A^{1}(t,z)\!-\!Q_p\!=\!0,\\
\rho\partial_t u^{1}(t,z)+\rho u^0(t,z)\partial_z u^{1}(t,z)+\beta u^{1}(t,z)+\partial_z P^0(t,z)=0,  \\
R\partial_t P^{1}(t,z)-KP^{1}(t,z)^2-KP^{1}(t,z)\left(RQ_p+\widetilde{P}\right)=0,
\end{cases}
\end{eqnarray}
\newpage
\begin{flushleft}
\textbf{As A2}
\end{flushleft}
\begin{eqnarray}\label{br10bbis}
\begin{cases}
\partial_t\eta^{1}(t,z)=\zeta^0(t,z),\\
\widehat{\alpha}A^{0}(t,z)\partial_{t}\zeta^{1}(t,z)+\widetilde{k}\zeta^{1}(t,z)+\kappa\eta^0(t,z)-A^0(t,z)P^0(t,z)+A^0(t,z)\widetilde{P}=0, \\
\left(\widetilde{h}+a(t)+\eta^0(t,z)\right)\partial_t A^{1}(t,z)+\left(a'(t)+\zeta^0(t,z)\right)A^{1}(t,z)-Q_p\\
\ \ \ \ \ \ \ \ \ \ \ \ \ \ \ \ \ \ \ \ \ \ \ \ \ \ \ \ \ \ \ \ \ \ \ \ \ \ \ \ \ \ \ \ \ \ \ \ \ \ \ \ \ \ \ \ \ \ \ \ \ \ \ \ \ \ \ +\displaystyle{\frac{1}{R}}\left(P^0(t,z)-\widetilde{P}\right)=0, \\
\rho\partial_t u^{1}(t,z)+\rho u^0(t,z)\partial_z u^{1}(t,z)+\beta u^{1}(t,z)+\partial_z P^0(t,z)=0,  \\
R\partial_t P^{1}(t,z)-KP^{1}(t,z)^2-KP^{1}(t,z)\left(RQ_p+\widetilde{P}\right)=0,
\end{cases}
\end{eqnarray}
with the following initial conditions
\begin{align}\label{br11bbis}
u(0,z)^{1}=f(z), \ \ & \eta(0,z)^{1}=g(z), \ \ \zeta(0,z)^{1}=q(z),\notag \\
P(0,z)^{1}&=b(z), \ \ A(0,z)^{1}=h(z).
\end{align}
Henceforth the analysis of the approximating systems is very similar, then we will specify the model only in case of  significant differences.

We want to show that there exist approximating solutions to the As A1 \eqref{br9bbis} and As A2 \eqref{br10bbis}.\\ First of all we observe that the solution of equation $\eqref{br9bbis}_5$ is defined by
\begin{equation}\label{br12bbis}
P^1(t,z)=\displaystyle{\frac{b(z)e^{\mathcal{C}t}}{1+\displaystyle{\frac{K}{R\mathcal{C}}}b(z)\left(1-e^{\mathcal{C}t}\right)}}.
\end{equation}
Then, by using the regularity conditions \eqref{br11bbis} we get
\begin{equation}\label{br13bbis}
P^1(t,z)\in C^{\infty}((0,T_0); \mathcal{H}^{s}([0,L]))
\end{equation}
and
\begin{equation}\label{br14bbis}
\partial_z P^0(t,z)\in C^{\infty}((0,T_0); \mathcal{H}^{s-1}([0,L])).
\end{equation}
Hence, by considering that the equation $\eqref{br9bbis}_4$ with initial condition $\eqref{br11bbis}_1$ defines a first order hyperbolic system, we can apply the Theorem \ref{iperb} with $M=[0,L]$ and $\mathbf{f}(x,t)=\partial_zP^1(t,z)$. Since every hypothesis of the Theorem is satisfied, we can conclude that there exists a unique solution $u^1$ such that 
\begin{equation}\label{br15bbis}
u^1(t,z)\in C([0,T_0); \mathcal{H}^{s-1}([0,L]))\cap C^1((0,T_0); \mathcal{H}^{s-2}([0,L])).
\end{equation}
Moreover we know that equations $\eqref{br9bbis}_3$ and $\eqref{br10bbis}_3$ with initial condition $\eqref{br11bbis}_5$ admit the following solution
\begin{equation}\label{br16bbis}
A^1(t,z)=h(z)e^{-G^0(t,z)t}+e^{-G^0(t,z)t}\int_0^t H^0(s,z)e^{G^0(s,z)s}\,ds,
\end{equation}
where
\begin{equation}\label{br17bbis}
G^0(t,z)=\displaystyle{\frac{\displaystyle{a'(t)+q(z)+f(z)}}{\widetilde{h}+a(t)+g(z)}}, \ \ \ \ H^0(t,z)=\displaystyle{\frac{Q_p}{\widetilde{h}+a(t)+g(z)}},
\end{equation}
for the system \eqref{br1bbis} and
\begin{equation}\label{br18bbis}
G^0(t,z)=\displaystyle{\frac{\displaystyle{a'(t)+q(z)+f(z)}}{\widetilde{h}+a(t)+g(z)}}, \ \ \ \ H^0(t,z)=\displaystyle{\frac{RQ_p-b(z)-\widetilde{P}}{R(\widetilde{h}+a(t)+g(z))}},
\end{equation}
for the system \eqref{br2bbis}. By using the assumptions \eqref{br11bbis} and the fact that $a(t)\in C^{\infty}([0,T_0])$, we obtain 
\begin{equation}\label{br19bbis}
A^1(t,z)\in C^{\infty}((0,T_0); \mathcal{H}^{s}([0,L])).
\end{equation}
Likewise, we show that the solutions to the initial value problems $\eqref{br9bbis}_2,$ $\eqref{br11bbis}_3$ and $\eqref{br9bbis}_1, \eqref{br11bbis}_2$ are respectively
\begin{align}\label{br20bbis}
\zeta^1(t,z)&=q(z)e^{-F^0(t,z)t}+\frac{1}{\widehat{\alpha}}e^{-F^0(t,z)t}\int_0^t\left(J^0(s,z)\right)e^{F^0(s,z)s}\,ds, \notag\\
\eta^1(t,z)&=g(z)+q(z)t,
\end{align}
where 
\begin{equation}\label{BR20bbis}
F^0(z,t)=\displaystyle{\frac{\widetilde{k}}{\widehat{\alpha}}h(z)}, \ \ \ \ J^0(z,t)=b(z)-\displaystyle{\frac{\kappa g(z)}{h(z)}}-\widetilde{P}.
\end{equation}
By the regularity of the initial conditions \eqref{br11bbis} we can conclude that
\begin{equation}\label{br21bbis}
\zeta^1(t,z)\in C^{\infty}((0,T_0); \mathcal{H}^{s}([0,L])),\ \ \eta^1(t,z)\in C^{\infty}((0,T_0); \mathcal{H}^{s}([0,L])).
\end{equation}
Therefore, the previous statements allow us to infer that there exist a unique solution
\begin{equation}\label{br22bbis}
\mathcal{X}^1(t,z)=(u^1(t,z), \eta^1(t,z), \zeta^1(t,z), P^1(t,z), A^1(t,z)),
\end{equation}
to the system \eqref{br9bbis} (\eqref{br10bbis} respectively) such that 
\begin{align}\label{br23bbis}
u^1(t,z)\in C([0,T_0); \mathcal{H}^{s-1}([0,L]))&\cap C^1((0,T_0); \mathcal{H}^{s-2}([0,L])),\notag \\
\eta^1(t,z)\in C^{\infty}((0,T_0); \mathcal{H}^{s}([0,L])), \ \ & \zeta^1(t,z)\in C^{\infty}((0,T_0); \mathcal{H}^{s}([0,L])), \notag \\
P^1(t,z)\in C^{\infty}((0,T_0); \mathcal{H}^{s}([0,L])), \ \ & A^1(t,z)\in C^{\infty}((0,T_0); \mathcal{H}^{s}([0,L])).
\end{align}
\begin{flushleft}
\textbf{$\mathbf{(n+1)}$-th STEP.}
\end{flushleft}
Now we want to prove there exists a unique solution for the systems \eqref{br1bbis} and \eqref{br2bbis}. In order to do that we recall that the \textit{n}-th iteration guarantees that there exists a unique solution $\mathcal{X}^n(t,z)$ to the approximating systems such that
\begin{align}\label{br30bbis}
u^n(t,z)\in C([0,T_0); \mathcal{H}^{s-1}([0,L]))&\cap C^1((0,T_0); \mathcal{H}^{s-2}([0,L])),\notag \\
\eta^n(t,z)\in C^{\infty}((0,T_0); \mathcal{H}^{s}([0,L])), \ \ & \zeta^n(t,z)\in C^{\infty}((0,T_0); \mathcal{H}^{s}([0,L])), \notag \\
P^n(t,z)\in C^{\infty}((0,T_0); \mathcal{H}^{s}([0,L])), \ \ & A^n(t,z)\in C^{\infty}((0,T_0); \mathcal{H}^{s}([0,L])).
\end{align}
Moreover, as well as the basic step, we know that 
\begin{equation}\label{br26bbis}
P^{n+1}(t,z)=\displaystyle{\frac{b(z)e^{\mathcal{C}t}}{1+\displaystyle{\frac{K}{R\mathcal{C}}}b(z)\left(1-e^{\mathcal{C}t}\right)}},
\end{equation}
then, since at every iteration the solution of the Riccati equation depends only on the initial data, we get
\begin{equation}\label{br27bbis}
P^{n+1}(t,z)\in C^{\infty}((0,T_0); \mathcal{H}^{s}([0,L]))
\end{equation}
and 
\begin{equation}\label{br28bbis}
\partial_z P^n(t,z)\in C^{\infty}((0,T_0); \mathcal{H}^{s-1}([0,L])).
\end{equation}
By considering \eqref{br28bbis}, we observe that the equation $\eqref{br1bbis}_4$ fulfills again the hypothesis of the Theorem \ref{iperb} which implies that there exists a unique solution 
\begin{equation}\label{br29bbis}
u^{n+1}(t,z)\in C([0,T_0); \mathcal{H}^{s-1}([0,L]))\cap C^1((0,T_0); \mathcal{H}^{s-2}([0,L])).
\end{equation}
This allow us to conclude that 
\begin{align}\label{br34bbis}
\zeta^{n+1}(t,z)&\in C([0,T_0); \mathcal{H}^{s-1}([0,L]))\cap C^1((0,T_0); \mathcal{H}^{s-2}([0,L])),\notag \\
\eta^{n+1}(t,z)&\in C([0,T_0); \mathcal{H}^{s-1}([0,L]))\cap C^1((0,T_0); \mathcal{H}^{s-2}([0,L])),\notag \\
A^{n+1}(t,z)&\in C([0,T_0); \mathcal{H}^{s-1}([0,L]))\cap C^1((0,T_0); \mathcal{H}^{s-2}([0,L])),
\end{align}
where
\begin{align}\label{br31bbis}
\zeta^{n+1}(t,z)&=q(z)e^{-F^n(t,z)t}+\frac{1}{\widehat{\alpha}}e^{-F^n(t,z)t}\int_0^t\left(J^n(s,z)\right)e^{F^n(s,z)s}\,ds,\notag \\
\eta^{n+1}(t,z)&=g(z)+\int_0^t\,\zeta^n(s,z)\,ds, 
\end{align}
with 
\begin{equation}\label{BR31bbis}
F^n(z,t)=\displaystyle{\frac{\widetilde{k}}{\widehat{\alpha}}A^n(t,z)}, \ \ J^n(z,t)=P^n(t,z)-\displaystyle{\frac{\kappa \eta^n(t,z)}{A^n(t,z)}}-\widetilde{P},
\end{equation}
and 
\begin{equation}\label{BR32bbis}
A^{n+1}(t,z)=h(z)e^{-G^n(t,z)t}+e^{-G^n(t,z)t}\int_0^t H^n(s,z)e^{G^n(s,z)s}\,ds,
\end{equation}
with
\begin{equation}\label{br32bbis}
G^n(t,z)=\displaystyle{\frac{\displaystyle{a'(t)+\zeta^n(t,z)+u^n(t,z)}}{\widetilde{h}+a(t)+\eta^n(t,z)}},\ \ H^n(t,z)=\displaystyle{\frac{Q_p}{\widetilde{h}+a(t)+\eta^n(t,z)}},
\end{equation}
in the approximating system A1 and
\begin{equation}\label{br33bbis}
G^n(t,z)=\displaystyle{\frac{\displaystyle{a'(t)+\zeta^n(t,z)+u^n(t,z)}}{\widetilde{h}+a(t)+\eta^n(t,z)}},\ \ H^n(t,z)=\displaystyle{\frac{RQ_p-P^n(t,z)-\widetilde{P}}{R(\widetilde{h}+a(t)+\eta^n(t,z))}},
\end{equation}
in the approximating system A2.\\
Now we are ready to outline the main result of this section in the following proposition.
\begin{proposition}[\textbf{Existence of the approximating solution}]\label{proposition1}
There exists a unique solution
\begin{equation}
\mathcal{X}=(u^{n+1}(t,z), \eta^{n+1}(t,z), \zeta^{n+1}(t,z), P^{n+1}(t,z), A^{n+1}(t,z)),
\end{equation}
to the problem \eqref{br1bbis}, \eqref{br3bbis} (\eqref{br2bbis}, \eqref{br3bbis}), where the sequences of functions are such that 
\begin{equation}\label{br34bbis}
P^{n+1}(t,z)\in C^{\infty}((0,T_0); \mathcal{H}^{s}([0,L])),
\end{equation}
and $u^{n+1}(t,z)$, $\eta^{n+1}(t,z)$, $\zeta^{n+1}(t,z)$, $A^{n+1}(t,z)$ are in
\begin{equation}
C([0,T_0); \mathcal{H}^{s-1}([0,L]))\cap C^1((0,T_0); \mathcal{H}^{s-2}([0,L])),
\end{equation}
for every $n$.
\end{proposition}

\subsection{Convergence of the iterative schemes}

In this section we focus on the second part of the proof of the Theorems \ref{TeoremaBr1} and \ref{TeoremaBr2}. From the Proposition \ref{proposition1} we know that there exists a unique solution for the approximating systems associated to Model A1 and Model A2 respectively and now we want to prove the strong convergence of these solutions to classical solutions of the systems \eqref{br1bis} and \eqref{br2bis}. \\
Due to its nonlinear behaviour, the equation $\eqref{br1bis}_4$ needs to be treated carefully. To this aim we recall the comprehensive analysis developed in \cite{donmaro} where the convergence of the approximating sequence $u^n(t,z)$ is proved by using high order energy estimates in the space $\mathcal{H}^s, s>7/2$. \\
The convergence of the other sequences of functions will be discuss in terms of equicontinuity.

\subsubsection{Analysis and convergence of the sequence $\pmb{\{u^n\}}$}\label{convU}

As we have already stated before, the sequence $u^n(t,z)$ is widely studied in \cite{donmaro} where the CSF velocity equation is exactly the same that appears in both Models A1 and A2, i.e.
\begin{equation}\label{br33}
\rho\partial_t u(t,z)+\rho u(t,z)\partial_z u(t,z)+\beta u(t,z)+\partial_z P(t,z)=0, 
\end{equation}
with initial condition
\begin{equation}\label{br33bis}
u(0,z)=u_0(z)=f(z)\in\mathcal{H}^s([0,L]),
\end{equation}
with $s>\displaystyle{\frac{5}{2}}$.
Then, based on the results obtained in \cite{donmaro}, we define the energy in order to perform a priori estimates. We proceed by setting the following notation
\begin{align}\label{br34}
u=u^{n}, \ \ v&=u^{n+1}-u_0, \ \  P=P^{n}, \notag \\
v_{\alpha}&=\frac{\partial^{\alpha}}{\partial z^{\alpha}}v, 
\end{align}
for any $\alpha<s$. For the sake of simplicity we will not indicate the space and time dependence in this section.\\
Now we define 
\begin{equation}\label{br36}
E^{\alpha}(t):=\frac{1}{2}\int_{0}^L\,\rho v_{\alpha}v_{\alpha}\, dz
\end{equation}
and we apply the operator $\displaystyle{\frac{\partial^{\alpha}}{\partial z^{\alpha}}}$ to the equation $\eqref{br33}$ and we get 
\begin{equation}\label{br37}
\rho\partial_t v_{\alpha}\!+\!\rho u\partial_z v_{\alpha}\!+\!\beta v_{\alpha}\!=\! \frac{\partial^{\alpha}}{\partial z^{\alpha}}(-\partial_z P\!-\!\rho u\partial_z f(z)\!-\!\beta f(z))\!-\!\rho\left(\frac{\partial^{\alpha}}{\partial z^{\alpha}}u\right)\partial_z v,
\end{equation}
\begin{equation}\label{br38}
v_{\alpha}(0,z)=0,
\end{equation} 
and 
\begin{align}\label{br39}
\frac{dE^{\alpha}}{dt}(t):=\int_{0}^L\,\rho \partial_t v_{\alpha}v_{\alpha}&\, dz= -\int_{0}^L\, \bigg[\frac{\partial^{\alpha}}{\partial z^{\alpha}}(\partial_z P+\rho u\partial_z f(z)+\beta f(z))\notag \\
&+\rho u\partial_z v_{\alpha}+\beta v_{\alpha}+\rho\left(\frac{\partial^{\alpha}}{\partial z^{\alpha}}u\right)\partial_z v\bigg] v_{\alpha}\,dz.
\end{align} 
By the detailed estimates in \cite{donmaro} we know that 
\begin{equation}\label{br40}
\frac{d\,E^{\alpha}(t)}{d\,t}\leq C_1\left(1+\left\|v_{\alpha}\right\|^2_2\right),
\end{equation}
where $C_1=C_1\left(\left\|u\right\|^2_{L^{\infty}_t\mathcal{H}^{s-1}_z}\right)$ is a locally bounded function of its argument.\\
Summing up \eqref{br40} over $\alpha<s-1$ and integrating over $(0,t)$, where $t\in(0,T)$, $T\leq T_0$, we find
\begin{align}\label{br41}
\sum_{\alpha<s-1}E^{\alpha}(t)\leq C_1(1+\left\|v\right\|_{L^{\infty}_t\mathcal{H}^{s-1}_z}^2)T.
\end{align}
Taking the $\sup$ over $t\in(0,T)$ in \eqref{br41} we get
\begin{equation}\label{br42}
\left\|v\right\|^2_{L^{\infty}_t\mathcal{H}^{s-1}_z}\leq C_1(\left\|u\right\|_{L^{\infty}_t\mathcal{H}^{s-1}_z})(1+\left\|v\right\|^2_{L^{\infty}_t\mathcal{H}^{s-1}_z})T,
\end{equation}
then
\begin{equation}\label{br43}
\left\|v\right\|^2_{L^{\infty}_t\mathcal{H}^{s-1}_z}\leq\frac{C_1T}{1-C_1T}:=C_2^2T.
\end{equation}
By using the notation \eqref{br34} this means that
\begin{equation}\label{br44}
\left\|u^{n+1}-u_0\right\|_{L^{\infty}_t\mathcal{H}^{s-1}_z}=\left\|u^{n+1}-f(z)\right\|_{L^{\infty}_t\mathcal{H}^{s-1}_z}\leq C_2(\left\|u^{n}\right\|_{L^{\infty}_t\mathcal{H}^{s-1}_z})\sqrt{T}.
\end{equation}
We take
\begin{equation}\label{br45}
r_0>\left\| f(z)\right\|_{\mathcal{H}^s_z}+C_2(\left\| f(z)\right\|_{\mathcal{H}^s_z}),
\end{equation}
and
\begin{equation}\label{br46}
\sqrt{T}<\left(\sup_{0\leq r\leq r_0}C_2(r)\right)^{-1}(r_0-\left\| f(z)\right\|_{\mathcal{H}^s_z}).
\end{equation}
Moreover we can show that for every $n\geq 0$ 
\begin{equation}\label{br47}
\left\| u^{n+1}\right\|_{L^{\infty}_t\mathcal{H}^{s-1}_z}\leq \left\| f(z)\right\|_{\mathcal{H}^{s}_z}+\sup_{0\leq r\leq r_0} C_2(r)\sqrt{T}<r_0,
\end{equation}
and
\begin{align}\label{br48}
\left\|\partial_t u^{n}\right\|_{L^{\infty}_t\mathcal{H}^{s-1}_z}&\leq \left(\left\| f(z)\right\|_{\mathcal{H}^{s}_z}+\sup_{0\leq r\leq r_0} C_2(r)\sqrt{T}\right)^2+\mathit{\overline{C}}\|b(z)\|_{\mathcal{H}^s_z}\notag\\
&+\left\lvert\frac{\beta}{\rho}\right\rvert\left(\left\| f(z)\right\|_{\mathcal{H}^{s}_z}+\sup_{0\leq r\leq r_0} C_2(r)\sqrt{T}\right),
\end{align}
where 
\begin{equation}\label{cost}
\mathit{\overline{C}}=\mathit{\overline{C}}\left(\mathcal{C}, T\right).
\end{equation}

In order to prove the convergence of the sequence $\{u^n\}_{n=0}^{\infty}$, we subtract two subsequent equations $\eqref{br1bbis}_4$ and we obtain
\begin{equation}\label{br49}
\rho\partial_t(u^{n+1}-u^n)+\rho u^n\partial_z(u^{n+1}-u^n)+\beta(u^{n+1}-u^n)= G_n,
\end{equation}
where
\begin{equation}\label{br51}
G_n:=\rho(u^{n-1}-u^{n})\partial_z u^n.
\end{equation}
By carrying out estimates similar to the previous ones, we obtain
\begin{align}\label{br52}
\left\| u^{n+1}-u^n\right\|_{L^{\infty}_t L^2_z}\leq C_3(\left\| f(z)\right\|_{\mathcal{H}^s_z},r_0,T)\sqrt{T}\left\| u^n-u^{n-1}\right\|_{L^{\infty}_t \mathcal{H}^{s-2}_z}
\end{align}
and we find that
\begin{equation}\label{br53}
\sum_{n=1}^{\infty} \left\| u^{n+1}-u^n\right\|_{L^{\infty}_t L^2_z}\leq\sum_{n=0}^{\infty}a^n\left\| u^{1}-u^0\right\|_{L^{\infty}_t \mathcal{H}^{s-2}_z}<\infty,
\end{equation}
where we assumed $T$ so small that
\begin{equation}\label{br54}
a:= C_3(\left\| f(z)\right\|_{\mathcal{H}^s_z},r_0,T)\sqrt{T}<1.
\end{equation}
The convergence \eqref{br53} implies that there exists $u\in L^{\infty}\left((0,T); L^2([0,L])\right)$ such that
\begin{equation}\label{br56}
u^n\rightarrow u \ \ \mbox{strongly in} \ \  L^{\infty}((0,T); L^2([0,L]))  \ \ \mbox{for some} \ \  u.
\end{equation}
Therefore, by \eqref{br47}, \eqref{br56} and by applying the interpolation inequality \eqref{ar2}, for $0<s'<s-1$, we obtain the following strong convergence
\begin{equation}\label{br57}
u^n\rightarrow u  \ \ \mbox{strongly in} \ \ C([0,T]; \mathcal{H}^{s'}([0,L])),
\end{equation}
and by the Sobolev embedding theorem we get
\begin{equation}\label{br58}
u^n\rightarrow u \ \ \mbox{strongly in}\ \ C([0,T]; C^1([0,L])).
\end{equation}
Moreover, from \cite{donmaro} (see Section 5.2) we are able to conclude that 
\begin{equation}\label{br59}
\partial_t u^n\rightarrow\partial_t u \ \ \mbox{strongly in} \ \ C((0,T]); C([0,L])).
\end{equation}
Because of the strong convergence \eqref{br58} and \eqref{br59}, we can assert that $u\in C^1((0,T]\times[0,L])$ is a classical solution of the equation $\eqref{br1bis}_4$.

\subsubsection{Convergence of the sequences $\pmb{\{P^n\}}, \pmb{\{\eta^n\}}, \pmb{\{\zeta^n\}}$ and $\pmb{\{A^n\}}$}

In this section, we shall show the convergence of $\{P^n\}, \{\eta^n\}, \{\zeta^n\}$ and $\{A^n\}$, and due to the different structure of the equations they satisfy we cannot employ the same methods as for $u^n$. Therefore we need to resort to other methods for the convergence proof, in particular we recall the following Ascoli's Theorem (\cite{stampacchia}).
\begin{theorem}[\textbf{Ascoli's Theorem}]\label{ascoli}
Let $X$ be a compact metric space, $Y$ a Banach space, and $C(X,Y)$ the Banach space of continuous functions from $X$ to $Y$ with the sup norm. A subset $\Omega$ of $C(X,Y)$ has compact closure in $C(X,Y)$ if and only if 
\begin{enumerate}[(a)]
\item $\Omega$ is equicontinuous, and
\item for every $x$ in $X$, the set $\Omega(x)=\{f(x): f\in\Omega\}$ has compact closure in $Y$.
\end{enumerate}
\end{theorem}
We want to apply the Theorem \ref{ascoli} to our sequences, more precisely we have to show that every sequence of functions is equicontinuous (see \eqref{equic2}). By using the Proposition \ref{proposition1} we obtain that
\begin{align}\label{br60}
\| P^{n+1}(t,z)-P^{n+1}(s,z)\|_{L^{\infty}_t L^2_z}&\leq\|\int_s^t\,\partial_{\tau}P^{n+1}(\tau,z)\,d\tau\|_{L^{\infty}_t \mathcal{H}^{s-1}_z}\notag \\
&\leq\lvert t-s\rvert\|\partial_{\tau}P^{n+1}(\tau,z)\|_{L^{\infty}_t \mathcal{H}^{s-1}_z}, \notag \\
\|\eta^{n+1}(t,z)-\eta^{n+1}(s,z)\|_{L^{\infty}_t L^2_z}&\leq\|\int_s^t\,\partial_{\tau}\eta^{n+1}(\tau,z)\,d\tau\|_{L^{\infty}_t \mathcal{H}^{s-1}_z}\notag \\
&\leq\lvert t-s\rvert\|\partial_{\tau}\eta^{n+1}(\tau,z)\|_{L^{\infty}_t \mathcal{H}^{s-1}_z}, \notag \\
\|\zeta^{n+1}(t,z)-\zeta^{n+1}(s,z)\|_{L^{\infty}_t L^2_z}&\leq\|\int_s^t\,\partial_{\tau}\zeta^{n+1}(\tau,z)\,d\tau\|_{L^{\infty}_t \mathcal{H}^{s-1}_z}\notag \\
&\leq\lvert t-s\rvert\|\partial_{\tau}\zeta^{n+1}(\tau,z)\|_{L^{\infty}_t \mathcal{H}^{s-1}_z}, \notag \\
\|A^{n+1}(t,z)-A^{n+1}(s,z)\|_{L^{\infty}_t L^2_z}&\leq\|\int_s^t\,\partial_{\tau}A^{n+1}(\tau,z)\,d\tau\|_{L^{\infty}_t \mathcal{H}^{s-1}_z}\notag \\ 
&\leq\lvert t-s\rvert\|\partial_{\tau}A^{n+1}(\tau,z)\|_{L^{\infty}_t \mathcal{H}^{s-1}_z},
\end{align}
where $t, s \in (0, T)$ and $s<\tau<t$.\\
By taking into account \eqref{br26bbis} for the pressure estimate and by obtaining directly $\partial_\tau \eta, \partial_\tau \zeta$ and $\partial_\tau A$ from the equations $\eqref{br1bbis}_1, \eqref{br1bbis}_2$ and $\eqref{br1bbis}_3$ ($\eqref{br2bbis}_3$ for As A2) respectively, we can write
 \begin{align}\label{br61}
\|\partial_{\tau}P^{n+1}(\tau,z)\|_{L^{\infty}_t \mathcal{H}^{s-1}_z}&\leq \widehat{M}\left\lvert\frac{K}{R\mathcal{C}}\right\rvert, \notag \\
\|\partial_{\tau}\eta^{n+1}(\tau,z)\|_{L^{\infty}_t \mathcal{H}^{s-1}_z}&\leq \|\zeta^{n}(\tau,z)\|_{L^{\infty}_t \mathcal{H}^{s}_z}, \notag\\
\|\partial_{\tau}\zeta^{n+1}(\tau,z)\|_{L^{\infty}_t \mathcal{H}^{s-1}_z}&\leq \left\lvert\frac{1}{\widehat{\alpha}}\right\rvert\bigg[\lvert\kappa\rvert\|\eta^n(\tau,z)\|_{L^{\infty}_t \mathcal{H}^{s}_z}\left(\widetilde{M}+1\right)\notag\\
&+\lvert\widetilde{k}\rvert\|A^n(\tau,z)\|_{L^{\infty}_t \mathcal{H}^{s}_z}\left(\widetilde{M}
+\|P^n(\tau,z)\|_{L^{\infty}_t \mathcal{H}^{s}_z}+\left\lvert\widetilde{P}\right\rvert\right)\bigg], \notag\\
\|\partial_{\tau}A^{n+1}(\tau,z)\|_{L^{\infty}_t \mathcal{H}^{s-1}_z}&\leq \overline{M}\|G^n(\tau,z)\|_{L^{\infty}_t \mathcal{H}^{s}_z}\left(\|h(z)\|_{\mathcal{H}^s_z}+T\|H^n(\tau,z)\|_{L^{\infty}_t \mathcal{H}^{s}_z}\right)\notag\\
&+\|H^n(\tau,z)\|_{L^{\infty}_t \mathcal{H}^{s}_z},
\end{align}
where $\widehat{M}=\widehat{M}\left(\|b(z)\|_{\mathcal{H}^s_z}, T\right)$, $\widetilde{M}=\widetilde{M}\left(\|q(z)\|_{\mathcal{H}^s_z}, \|b(z)\|_{\mathcal{H}^s_z}, T\right)$ and $\overline{M}=\overline{M}\left(T, r_0, \|g(z)\|_{\mathcal{H}^s_z}\right)$. \\
By considering the Proposition \ref{proposition1}, we observe that
\begin{align}\label{br62}
\|\partial_{\tau}P^{n+1}(\tau,z)\|_{L^{\infty}_t \mathcal{H}^{s-1}_z}&\leq \widehat{M'}, \ \ \|\partial_{\tau}\eta^{n+1}(\tau,z)\|_{L^{\infty}_t \mathcal{H}^{s-1}_z}\leq \mathcal{M'}, \notag \\
\|\partial_{\tau}\zeta^{n+1}(\tau,z)\|_{L^{\infty}_t \mathcal{H}^{s-1}_z}&\leq \widetilde{M'}, \ \ \|\partial_{\tau}A^{n+1}(\tau,z)\|_{L^{\infty}_t \mathcal{H}^{s-1}_z}\leq \overline{M'}.
\end{align}
If we denote with $M=\min\left\{\widehat{M'}, \mathcal{M'}, \widetilde{M'}, \overline{M'}\right\}$ and we plug properly \eqref{br62} in \eqref{br60}, we conclude that 
\begin{align}\label{br63}
\| P^{n+1}(t,z)-P^{n+1}(s,z)\|_{L^{\infty}_t L^2_z}&\leq M \lvert t-s\rvert, \notag \\
\|\eta^{n+1}(t,z)-\eta^{n+1}(s,z)\|_{L^{\infty}_t L^2_z}&\leq M \lvert t-s\rvert, \notag \\
\|\zeta^{n+1}(t,z)-\zeta^{n+1}(s,z)\|_{L^{\infty}_t L^2_z}&\leq M \lvert t-s\rvert, \notag \\
\|A^{n+1}(t,z)-A^{n+1}(s,z)\|_{L^{\infty}_t L^2_z}&\leq M \lvert t-s\rvert.
\end{align}
This last result suggests the sequences $\{P^n\}, \{\eta^n\}, \{\zeta^n\}$ and $\{A^n\}$ are therefore equicontinuous (see Definition \ref{equic2}); hence is trivial to prove they are uniformly bounded, as well. \\
We can apply now the Theorem \ref{ascoli} which implies that the set 
\begin{displaymath}
\Psi(t,z)=\left\{P^n(t,z), \eta^n(t,z), \zeta^n(t,z), A^n(t,z); t\in[0,T], z\in[0,L]\right\}
\end{displaymath}
is relatively compact in 
\begin{displaymath}
C(X,Y)=C([0,T]; \mathcal{H}^{s-1}([0,L]))\cap C^1((0,T]; \mathcal{H}^{s-2}([0,L])),
\end{displaymath}
where $T<T_0$. \\
Then, the relative compactness of $\Psi$ allows us to conclude that 
 \begin{align}\label{br64}
P^{n}\rightarrow P \ \  &\mbox{strongly in} \ \  {L^{\infty}((0,T); L^2([0,L]))}, \notag \\ 
\eta^{n}\rightarrow\eta \ \  &\mbox{strongly in} \ \  {L^{\infty}((0,T); L^2([0,L]))}, \notag \\
\zeta^{n}\rightarrow\zeta \ \  &\mbox{strongly in} \ \  {L^{\infty}((0,T); L^2([0,L]))}, \notag \\
\partial_t\zeta^{n}\rightarrow\partial_t\zeta \ \  &\mbox{strongly in} \ \  {L^{\infty}((0,T); L^2([0,L]))}, \notag \\
A^{n}\rightarrow A\ \  &\mbox{strongly in} \ \  {L^{\infty}((0,T); L^2([0,L]))},
 \end{align}
and finally, by applying the interpolation inequality \eqref{ar2}, we get
\begin{align}\label{br65}
P^{n}\rightarrow P \ \  &\mbox{strongly in} \ \  C((0,T]); C([0,L])), \notag \\ 
\eta^{n}\rightarrow\eta \ \  &\mbox{strongly in} \ \  C((0,T]); C([0,L])), \notag \\
\zeta^{n}\rightarrow\zeta \ \ & \mbox{strongly in} \ \  C((0,T]); C([0,L])), \notag \\
\partial_t\zeta^{n}\rightarrow\partial_t\zeta \ \  &\mbox{strongly in} \ \  C((0,T]); C([0,L])), \notag \\
A^{n}\rightarrow A\ \ & \mbox{strongly in} \ \  C((0,T]); C([0,L])).
 \end{align}

\subsection{Existence and uniqueness of the solution}\label{localsol}

In the previous sections we attained the convergence of the approximating sequences of Models A1 and A2. This means that we are allowed to pass into the limit in the systems \eqref{br1bbis} and \eqref{br2bbis} and we obtain that 
\begin{equation}\label{br66}
\mathcal{X}=\left(u, P, \eta, \zeta, A\right)\in C((0,T]\times[0,L])
\end{equation}
is a classical solution of \eqref{Tbr1} (\eqref{Tbr21}).\\
The last result to be shown is that the solution is unique. We assume there is another solution $\Psi=\left(\overline{u}, \overline{P}, \overline{\eta}, \overline{\zeta}, \overline{A}\right)$ of the system \eqref{Tbr1} (\eqref{Tbr21}) with initial condition \eqref{Tbr2}. \\
By denoting 
\begin{equation}\label{br67}
\widehat{u}=u-\overline{u}, \ \ \widehat{P}=P-\overline{P}, \ \ \widehat{\eta}=\eta-\overline{\eta}, \ \ \widehat{\zeta}=\zeta-\overline{\zeta}, \ \ \widehat{A}=A-\overline{A},
\end{equation}
and after analogous calculations as performed in Section \ref{convU} (see in \cite{donmaro}), we end up with the inequalities
\begin{align}\label{br68}
\left\|\widehat{u}(t)\right\|_2^2\leq\widetilde{C_1}\int_0^t\,\left\|\widehat{u}(s)\right\|_2^2\,ds,\ \ & \left\|\widehat{P}(t)\right\|_2^2\leq\widetilde{C_2}\int_0^t\,\left\|\widehat{P}(s)\right\|_2^2\,ds,\notag \\
\left\|\widehat{\eta}(t)\right\|_2^2\leq\widetilde{C_3}\int_0^t\,\left\|\widehat{\eta}(s)\right\|_2^2\,ds, \ \ & \left\|\widehat{\zeta}(t)\right\|_2^2\leq\widetilde{C_4}\int_0^t\,\left\|\widehat{\zeta}(s)\right\|_2^2\,ds, \notag \\
\left\|\widehat{A}(t)\right\|_2^2\leq\widetilde{C_5}\int_0^t\,&\left\|\widehat{A}(s)\right\|_2^2\,ds,
\end{align}
with $\widehat{u}(0)=0, \widehat{P}(0)=0, \widehat{\eta}(0)=0, \widehat{\zeta}(0)=0$ and $\widehat{A}(0)=0$.
This yields 
\begin{equation}\label{br69}
\mathcal{X}-\Psi\equiv 0
\end{equation}
and concludes the proof of the Theorems \ref{TeoremaBr1} and \ref{TeoremaBr2}.

\section{PROOF OF THE THEOREM \ref{TeoremaBr3}}

This section is devoted to the problem of the global existence of solutions for Models A1 and A2. This issue has been already investigated in \cite{donmaro} and for the simplified CSF analyzed model it has been proved existence and uniqueness of a global solution under some restriction conditions on the initial data. As well as in \cite{donmaro} the core of the problem here is the Burgers-like behaviour of the equation
\begin{equation}\label{br70}
\rho\partial_t u(t,z)+\rho u(t,z)\partial_z u(t,z)+\beta u(t,z)+\partial_z P(t,z)=0,
\end{equation}
which affects the existence time of the solutions $\eta, \zeta$ and $A$ for the other equations that appear in the systems A1 and A2. This peculiarity allows us to apply in the present paper the results obtained in \cite{donmaro} with some proper variations related to the structure of the new models. 

\subsection{Global existence of solutions}

In order to evaluate if it is possible to prove global existence of solutions we want to draw attention to the following points.
\begin{enumerate}[(a)]
\item From the Section \ref{localsol} we know that there exists a unique local solution $\mathcal{X}=\left(u, P, \eta, \zeta, A\right)\in C((0,T]\times[0,L])$ where 
\begin{equation}\label{br71}
\sqrt{T}<\left(\sup_{0\leq r\leq r_0}C_2(r)\right)^{-1}(r_0-\left\| f(z)\right\|_{\mathcal{H}^s_z}).
\end{equation}
\item In \cite{donmaro} the authors analyze first the homogeneous equation 
\begin{equation}\label{abr72}
\partial_t u+u\partial_z u+\frac{\beta}{\rho}u=0,
\end{equation}
with the initial condition 
\begin{equation}\label{br73}
u(0,z)=u_0(z)=f(z)\in\mathcal{H}^s([0,L]),
\end{equation}
with $s>7/2$.
By applying the method of characteristics, they obtain a Riccati Cauchy problem satisfied by $\omega=u_z$ along the characteristic curve, and by classical ODE method they are able to conclude that the solution of \eqref{abr72}, \eqref{br73} exists globally if and only if 
\begin{equation}\label{br74}
f'(z)\geq-\displaystyle{\frac{\beta}{\rho}}.
\end{equation}
This condition lays the foundations for the global existence proof of the solution to the equation \eqref{br70}.
\item According to the analysis developed in \cite{donmaro}, we apply the method of characteristics to \eqref{br70} and we denote by $\omega=u_z$, which satisfies along the characteristic curve $\Gamma_{\lambda}$ the following ordinary nonlinear differential Cauchy problem
\begin{eqnarray}\label{copia1}
\begin{cases}
\displaystyle{\omega'+\omega^2+\frac{\beta}{\rho}\omega+\frac{\beta}{\rho}P_{zz}=0}, \\
\omega(0)=f'(\lambda),
\end{cases}
\end{eqnarray}
where $\lambda=z$ for $t=0$.\\
Therefore, the focus of the attention shifts to the study of a Riccati nonhomogeneous equation for which we need to find a particular solution, $\overline{\omega}(t)$. This problem has been already treated in \cite{donmaro} and we know that the existence of $\overline{\omega}(t)$ is guaranteed if and only if 
\begin{equation}\label{br77c}
\lVert P_{zz}(t,z)\rVert_{L^2_{t,z}}\leq\infty.
\end{equation}
Moreover, since we are dealing with physiological models, we are interested in a real particular solution which can be obtained if and only if
\begin{equation}\label{br77}
\sup_{t,z}\,\lvert P_{zz}\rvert\leq\frac{\beta}{4\rho},
\end{equation}
(for more details see \cite{donmaro}, Section 7.2).
\end{enumerate}
\vspace{8pt}

By taking into account the previous points we proceed in two steps, as follows. 
\begin{flushleft}
\textit{STEP 1.}
\end{flushleft}
\vspace{4pt}

First of all we recall that the solution of the equation
\begin{equation}\label{br75}
R\partial_t P(t,z)-KP(t,z)^2-KP(t,z)\left(RQ_p+\widetilde{P}\right)=0, 
\end{equation}
with initial condition $P(0,z)=b(z)\in\mathcal{H}^s$, is the following 
\begin{equation}\label{br76}
P(t,z)=\displaystyle{\frac{b(z)e^{\mathcal{C}t}}{1+\displaystyle{\frac{K}{R\mathcal{C}}}b(z)\left(1-e^{\mathcal{C}t}\right)}}.
\end{equation} 
Therefore, we deduce higher regularity for the pressure than the one achieved in \cite{donmaro}. In particular, if $b(z)>0$, this is true provided $t\in[0, T)$, where $T<T_0<\widetilde{T}_0$ and $\widetilde{T}_0=\displaystyle{\frac{1}{\mathcal{C}}\ln\left(\frac{R\mathcal{C}}{Kb(z)}+1\right)}$ is the blow up time (see Remark \ref{remark1}). By the previous assumption, in this case we are sure that 
\begin{equation}\label{EPS}
1+\displaystyle{\frac{K}{R\mathcal{C}}}b(z)\left(1-e^{\mathcal{C}t}\right)\geq \varepsilon>0.
\end{equation} 
\par
We observe that the condition \eqref{br77c} is automatically satisfied by taking into account \eqref{br76} and the local existence Theorems \ref{TeoremaBr1} and \ref{TeoremaBr2}.\\
Nevertheless, as explained in point (c), we need real solutions which are guaranteed if and only if the condition \eqref{br77} is fulfilled. In order to show the latter relation, by considering the solution \eqref{br76}, we observe that 
\begin{equation}\label{brPr1}
\lVert P_{zz}(t,z)\rVert_{L^{\infty}_t \mathcal{H}^{s-2}_z}\leq \lVert P(t,z)\rVert_{L^{\infty}_t \mathcal{H}^{s}_z}\leq \widehat{\mathcal{C}}\|b(z)\|_{\mathcal{H}^{s}_z},
\end{equation}
where 
\begin{itemize}
\item $\widehat{\mathcal{C}}=e^{\mathcal{C}T}$ if $b(z)<0$,
\item $\widehat{\mathcal{C}}=\displaystyle{\frac{e^{\mathcal{C}T}}{\varepsilon}}$ if $b(z)>0$ (cf. \eqref{EPS} and Remark \ref{remark1}),
\end{itemize}
and $T$ is the local existence time defined in \eqref{br71}. Therefore, the condition \eqref{br77} can be satisfied if and only if
\begin{equation}\label{brPr2}
\|b(z)\|_{\mathcal{H}^{s}_z}\leq\frac{\beta}{4\widehat{\mathcal{C}}_1\rho},
\end{equation}
when the pressure initial datum is negative, and
\begin{equation}\label{br83}
\|b(z)\|_{\mathcal{H}^{s}_z}\leq\varepsilon\frac{\beta}{4\widehat{\mathcal{C}}_1\rho},
\end{equation}
if $b(z)$ is positive. The constant that appears in the previous conditions is obtained by considering \eqref{br71} and is defined as $\widehat{\mathcal{C}}_1=\widehat{\mathcal{C}}_1\left(\|f(z)\|_{\mathcal{H}^s_z}\right)$.\\
\begin{remark}
We want to point out that, while in \cite{donmaro} the condition \eqref{br77} is satisfied on the interval $[0, \widetilde{T}]$ with $\widetilde{T}<T$, where $T$ is the local existence time, if and only if
\begin{equation}\label{br78}
\left\lVert f(z)\right\rVert_{\mathcal{H}^s_z}\leq\frac{\beta}{\rho},
\end{equation}
in the present analysis of Models A1 and A2, we are not forced to further restrict the condition \eqref{br74} which guarantees the global existence of a solution for the homogeneous Riccati equation \eqref{abr72} and, therefore, it represents a first necessary assumption for the global existence of a solution for the problem \eqref{copia1}. The relation \eqref{br77} is fulfilled, here, by imposing the costraint \eqref{brPr2} or \eqref{br83} on the initial datum of the pressure.
\end{remark}

In conclusion, the Cauchy problem \eqref{copia1} admits the following unique solution
\begin{equation}\label{solRiccati}
\omega(t)=\overline{\omega}(t)+\frac{1}{y(t)},
\end{equation}
where
\begin{equation}\label{solRiccati2}
y(t)=\displaystyle{e^{\int_0^t\left(2\overline{\omega}(\tau)+\frac{\beta}{\rho}\right)\,d\tau}}\left\{\frac{2}{f'(z)}+\int_0^t e^{-\int_0^t\left(2\overline{\omega}(s)+\frac{\beta}{\rho}\right)\,ds}\,d\tau\right\},
\end{equation}
if and only if the conditions \eqref{br74} and \eqref{brPr2} (or \eqref{br83}) are satisfied. This solution is defined for any time $t\in[0,T^{\ast}]$, such that 
\begin{equation}\label{br84}
u(t,z) \in C\left([0, T^{\ast}]; \mathcal{H}^{s-1}([0, L])\right)\cap C^1\left([0, T^{\ast}); \mathcal{H}^{s-2}([0,L])\right),
\end{equation}
with 
\begin{equation}\label{br85}
T^{\ast}\leq \min\left\{T, \widetilde{T}_0\right\}.
\end{equation}

\begin{flushleft}
\textit{STEP 2.}
\end{flushleft}
\vspace{4pt}
By taking into account the solution \eqref{solRiccati}, the relation \eqref{brPr1} and the conditions \eqref{brPr2},  \eqref{br83}, we can therefore write that
\begin{align}
\left\lVert u_z\right\rVert_{L^{\infty}_t\mathcal{H}^{s-2}_z}&\leq 1+\left\| P_{zz}\right\|_{L^{\infty}_t\mathcal{H}^{s-2}_z}\notag\\
&\leq 1+\frac{\beta}{4\widehat{\mathcal{C}}_1\rho}=M_1,
\end{align}
if $b(z)<0$, and
\begin{equation}
\left\lVert u_z\right\rVert_{L^{\infty}_t\mathcal{H}^{s-2}_z}\leq 1+\varepsilon\frac{\beta}{4\widehat{\mathcal{C}}_1\rho}=M_2,
\end{equation}
if $b(z)>0$.\\
In order to complete the proof of the Theorem \ref{TeoremaBr3}, we have just to recall the following theorem (\cite{majda}).
\begin{theorem}[A Sharp Continuation Principle]\label{TeoremaG}
Let be $u_0\in\mathcal{H}^s$ for some $s>\frac{5}{2}$ and $T'>0$ some given time. Assume that for any interval of classical existence $[0, T^{\ast}]$, $T^{\ast}\leq T'$ for u(t) solution of the equation \eqref{br70}, the following a priori estimate is satisfied:
\begin{equation}\label{TeoremaG1}
\left\lVert u_z\right\rVert_{L^{\infty}_t}\leq M, \ \ \ \ \ \ \ \ 0\leq t\leq T^{\ast},
\end{equation}
where $M$ is a fixed constant independent of $T^{\ast}$. \\
Then the classical solution $u(t)$ exists on the interval $[0, T']$, with $u(t)$ in $C\left([0,T'], \mathcal{H}^{s-1}([0, L])\right)\cap C^1\left([0,T'], \mathcal{H}^{s-2}([0,L])\right)$. Furthermore, $u(t)$ satisfies the a priori estimate
\begin{equation}\label{TeoremaG2}
\left\| u\right\|_{L^{\infty}_t\mathcal{H}^{s-2}_z}\leq C e^{MT^{\ast}}\left\| f(z)\right\|_{\mathcal{H}^s_z},
\end{equation}
for $T^{\ast}$ with $0\leq T^{\ast}\leq T'$ and the constants $C, M$ depend on the space interval $[0, L]$ and on some physics quantities.
\end{theorem}
Since the hypothesis of the Theorem \ref{TeoremaG} are satisfied, we can conclude that there exists a global unique solution $\mathcal{X}(t,z)=(u, P, \eta, \zeta, A)(t,z)$ to the systems of equations \eqref{Tbr1} and \eqref{Tbr21}. 

\section*{Acknowledgements}
The authors would like to deeply thank the referees for the several accurate and helpful suggestions which have been essential to improve the paper.

\end{document}